\documentclass[12pt,twoside]{amsart}
\usepackage{amssymb}

\newtheorem{thm}{Theorem}[section]
\newtheorem{lemma}[thm]{Lemma}
\newtheorem{prop}[thm]{Proposition}
\newtheorem{corol}[thm]{Corollary}
\newtheorem{obs}[thm]{Remark}
\theoremstyle{definition}
\newtheorem{defi}[thm]{Definition}
\newtheorem{notation}[thm]{Notation}
\newtheorem{condition}[thm]{Condition}
\newcommand{\oo}{{\mathcal O}}
\newcommand{\cE}{{\mathcal E}}
\newcommand{\ii}{{\mathcal I}} 
\newcommand{\cI}{{\mathcal I}}
\newcommand{\cN}{{\mathcal N}}
\newcommand{\cF}{{\mathcal F}}
\newcommand {\pp}{\mathbb{P}}
\newcommand {\AAA}{\mathbb{A}}

\newcommand{\al}{\alpha}
\newcommand{\eeps}{\varepsilon}

\newcommand{\im}{\operatorname{im}}
\newcommand{\coker}{\mbox{coker}}
\newcommand{\lra}{\longrightarrow}

\newcommand{\isom}{\cong}
\newcommand{\ffi}{\varphi}
\newcommand\codim{\operatorname{codim}}
\newcommand\id{\operatorname{id}}
\newcommand{\s}{\; | \;}
\newcommand{\HR}{Hartshorne-Rao module}
\newcommand{\ER}{\operatorname{Ext}_R}
\newcommand{\Proj}{\operatorname{Proj}} 
\newcommand{\chr}{\operatorname{char}} 
\newcommand{\mif}{\mbox{if} ~} 
\newcommand{\fall}{\mbox{for all} ~}
\newcommand{\Grass}{\operatorname{Grass}}
\newcommand{\adj}{\operatorname{adj}}

\newcommand {\ZZ}{\mathbb{Z}}

\begin{document}
\title[The Hilbert scheme of degree two curves and certain ropes]{The Hilbert
scheme of degree two curves and certain ropes}
\author[Uwe Nagel, Roberto Notari, Maria Luisa Spreafico]{Uwe Nagel,
Roberto Notari, Maria Luisa Spreafico} 
\address{Department of Mathematics,
University of Kentucky, 715 Patterson Office Tower, 
Lexington, KY 40506-0027, USA}
\email{uwenagel@ms.uky.edu} 
\address{Dipartimento di Matematica, Politecnico di Torino, I-10129 Torino,
Italy}
\email{notari@polito.it} 
\address{Dipartimento di Matematica, Politecnico di Torino, I-10129 Torino,
Italy} 
\email{spreafico@polito.it}

% \date{\today} 

\begin{abstract} We study families of ropes of any codimension that
  are supported on 
lines.  In particular, this includes all non-reduced curves of degree
two. 
We construct suitable smooth parameter spaces and conclude
that all ropes of fixed degree and genus lie in the same component of
the corresponding Hilbert scheme. We show that this component is
generically smooth if the genus is small enough unless the
characteristic of the ground field is two and the curves under
consideration have degree two.  In this case the component is
non-reduced.  

\end{abstract}
\subjclass{Primary 14C05, 14H10, 14H45; Secondary 13D02 }

\thanks{The authors warmly thank INDAM-GNSAGA for partial support. The
  first author was also partially supported by a Faculty Summer
  Research Fellowship from the University of Kentucky.} 
%%%%%%%%%%%%%%%%%%%%%%%%%%%%%%%%%%%%%
\maketitle
\tableofcontents
%%%%%%%%%%%%%%%%%%%%%%%%%%%%%%%%%%%%%%%%%%%%5
\section{Introduction}

In this paper we study families of certain curves of small degree and
arbitrary codimension including all curves of degree two. By a curve we always
mean a closed locally Cohen-Macaulay projective subscheme of pure
dimension one. If the arithmetic genus of a degree two curve  is
small enough then it is non-reduced, supported on a line, and
contained  in the second infinitesimal neighborhood of the
line. Curves with these properties are called ropes. 
 We
construct smooth parameter spaces for ropes and use 
them to study Hilbert schemes that contain ropes.  

There are only few results about general properties of the Hilbert
scheme as its existence \cite{Grothen-Hilb}, its connectedness
\cite{Harts-Hilb}, or its radius \cite{Reeves-rad}. More is known
about Hilbert schemes parametrizing particular schemes that often have
small codimension (cf., e.g., \cite{Piene-Schles}, \cite{MD-Piene},
\cite{Ellingsrud}, \cite{MDP-Hilb-scheme}, \cite{Nollet},
\cite{Nollet-Schl}). But, for 
example, it is not even known if the locus of space curves inside its
Hilbert scheme is connected (cf., e.g., \cite{Hartsh-conn-Hilb},
\cite{Perrin}).  In this
note we also consider particular schemes, certain curves of low
degree, but we allow for any codimension.  

Curves of degree one, i.e.\ lines, are parametrized by
Grassmanians. Here, we study the next case(s). A curve of degree two
is either an irreducible conic, a pair of two lines, or a double
line. If such a curve is not planar its (arithmetic) genus is at
most $-1$. If it is equal to $-1$ then the Hilbert scheme is well
understood. Its general curve is a pair of skew lines. Somewhat
surprisingly,  if the genus is smaller then the characteristic of the
ground field enters the picture. Let $Hilb_{d, g}(\pp^n_K)$ denote the
Hilbert scheme of $1$-dimensional subschemes of $\pp^n_K$ having
degree $d$ and arithmetic genus $g$ where $K$ is an algebraically
closed field. Of course, $Hilb_{d, g}(\pp^n_K)$ does not contain only
curves. 
Part of our results on degree two curves  (cf.\ Theorem
\ref{thm-double-line}) is summarized in the following statement.  

\begin{thm} \label{thm-intro-deg-two} 
Let $g \leq -2,\ n \geq 3$ be integers. Then there is a unique
component, say $H_{2, g}(\pp^n_K)$, of $Hilb_{2, g}(\pp^n_K)$ that
parametrizes all curves in $Hilb_{2, g}(\pp^n_K)$. This component is
generically smooth if and only if either $\chr K \neq 2$ or $\chr K =
2$ and $g = -2$. If $\chr K = 2$ and $g \leq -3$ then $H_{2,
  g}(\pp^n_K)$ is non-reduced.  
\end{thm}

Even in the classical case of space curves, i.e.\ $n = 3$, this
results seems new in case of characteristic two.  
 
The general curve of $H_{2, g}(\pp^n_K)$ is a double line if $g \leq
-2$. Note that a double line is a rope of degree two. The result above
is a particular case of our effort to understand families of ropes
that are supported on a line. In general, it turns out that all ropes
in $\pp^n$ of fixed degree $d$ and (arithmetic) genus $g$ lie in the
same irreducible component of  $Hilb_{d, g}(\pp^n_K)$ (cf.\ Lemma
\ref{lem-unique-comp}). We denote this component by $H_{d,
  g}(\pp^n_K)$. If $d \geq 3$ then there are many more curves of
degree $d$ than just ropes. But we show, as the case of degree two
curves suggests, if the genus is small enough then the general curve of
$H_{d, g}(\pp^n_K)$ is a rope. Moreover, we have.  

\begin{thm} \label{thm-intro-ropes} 
Let $d, g, n$ be integers such that $3 \leq d \leq n-1$. 
If $g \leq \min \{-3(d-1), d-n\}$ then the component $ H_{d, g}(\pp^n_K)$ is
generically smooth of dimension $ (n-1)(d+1-g) - (d-1)^2$. 
\end{thm} 

For a more precise result we refer to Theorem \ref{thm-main-ropes}. It is independent of the characteristic of the ground field $K$.

A crucial ingredient of the results above is the computation of the global sections of the normal sheaf of a rope $C$.  However, we were not able to compute $h^0(C, \cN_C)$
by using general methods like vanishing theorems. Indeed, the dependence on the characteristic suggests that this might be impossible. Instead we use a very concrete approach based on our previous results in \cite{NNS}. It requires information on the minimal free resolution of the rope and its structure sheaf. 

Let $C \subset \pp^n$ be a rope of degree $n - k < n$ supported on the line $L$. Then $C$ is intimately related to two matrices $A, B$ with entries in the coordinate ring $S$ of $L$ that fit into an exact sequence 
$$
0 \to \oplus_{j=1}^k \oo_{L}(-\beta_j-1)
\stackrel{B}{\lra} \oo_{L}^{n-1}(-1)
\stackrel{A}{\lra} \oplus_{i=0}^{r-k} \oo_{L}(\al_i-1)
\to 0.
$$ 
Invoking a change of coordinates we may assume that the line $L$ is defined by the ideal $(x_0,\ldots,x_{n-2})$. Then, according to \cite{NNS}, Theorem 2.4, the homogeneous ideal of the rope is 
$$
I_C = ((I_L)^2,\; [x_0,\ldots,x_{n-2}] \cdot B). 
$$
Using this information we compute in Section \ref{sec-res} the graded minimal free resolution of $C$ including a description of the maps (Theorem \ref{Cres}). This is achieved by induction on the degree of the rope and based on a {\it non-minimal} resolution of the ideal $(I_L)^2$. 

Section \ref{sec-struc-sheaf} is devoted to determine the global
sections of the structure sheaf $\oo_C$ of the rope $C$. Using the
matrix $A$ associated to $C$ we construct explicitly global sections
that form a basis of $H^0_*(C, \oo_C) := \oplus_{j \in \ZZ} H^0(C,
\oo_C (j))$ as $S$-module (cf.\ Theorem \ref{thm-desc-gl-sec}). This
generalizes Migliore's result (\cite{mig}) for double lines in
$\pp^3$. It also serves as the basis for determining the module
structure of $H^0_*(C, \oo_C)$ over the coordinate ring $R =
K[x_0.\ldots,x_n]$ of $\pp^n$. Actually, we even compute the minimal
free resolution of $H^0_*(C, \oo_C)$ in Proposition
\ref{prop-res-struc}.  

Let 
$$
0 \to G_n \stackrel{d_n'}{\lra} \dots \lra G_2 \stackrel{d_2'}{\lra} G_1
\stackrel{d_1'}{\lra} I_C \to 0
$$
be the minimal free resolution of $C$. Then the global sections of its
normal sheaf can be computed by  
$$
H^0(C, {\mathcal N}_C) = [\ker(\mbox{Hom}(G_1, H^0_*(C, {\mathcal
O}_C)) \stackrel{(d_2')^{*}}{\lra} \mbox{Hom}(G_2, H^0_*(C,
{\mathcal O}_C)))]_0 
$$

This amounts to solving a system of linear equations over the
polynomial ring $R$ and is the content of Section
\ref{sec-normal-sheaf}. We are able to do this if $C$ has degree two,
i.e.\ for a double line. In this case, the solutions of the linear
system have a particular property that we formalize in Condition
\ref{cond}. For an arbitrary rope $C$ we compute all solutions of the
linear system that satisfy Condition \ref{cond}.  

In Section \ref{sec-para-spaces} we turn to families of ropes. Using the matrices $A, B$ associated to a rope $C$ we construct two types of smooth parameter spaces (cf.\ Propositions \ref{prop-dim-Val} and \ref{prop-Vbe}). Then a deformation argument allows to conclude that all ropes of fixed degree and genus lie in the same irreducible component $H_{d, g}(\pp^n)$ of the corresponding Hilbert scheme. Combined with the information of Section \ref{sec-normal-sheaf} on the tangent sheaf of $Hilb_{2, g}(\pp^n)$ at a point corresponding to a double line we obtain in particular Theorem \ref{thm-intro-deg-two}. 

For ropes of degree $d \geq 3$ the situation is more complicated because the dimension of the corresponding tangent space depends on the rope. In Proposition \ref{prop-small-coho} we determine the ropes with smallest cohomology. As expected by semicontinuity, it turns out that these ropes are indeed general in $H_{d, g}(\pp^n)$. Moreover, we show that the global sections of the normal sheaf of a general rope satisfy Condition \ref{cond}. This leads to Theorem \ref{thm-intro-ropes}. 
\smallskip

\begin{center}
{\bf Acknowledgement}
\end{center} During the preparation of this paper the first author
profited by a visit at the  Politecnico di Torino while the second and
third author profited by a stay at the University of Kentucky. We
thank these institutions for their hospitality. 

%%%%%%%%%%%%%%%%%%%%%%%%%%%%%%%%%%%%%%%%%%%%%%%%%%%%%%%%%%%%%%%%%%
\section{Characterizations of ropes on a line} \label{sec-char} 

In this section we introduce some notation and recall results from \cite{NNS} that we use later on. 

Throughout the paper a curve will always be a locally Cohen-Macaulay
$1$-dimensional projective subscheme.  A rope $C$ is a curve that is
supported on a smooth curve  $Y$  such that its homogeneous ideal $I_C$
contains $I_Y^2$. 
In this note we will only consider ropes whose support is a line. We often will refer to these curves just as ropes,  not making the assumption on the support explicit. A $d$-rope is just a rope of degree $d$.  We may and will assume that the supporting line $L$ is defined by the ideal $I_L = (x_0,\ldots,x_{n-2})$. In order to stress the particular role played by the line we denote the coordinate ring of $\pp^n$ by $R = K[x_0,\ldots,x_r, t, u]$ where $r := n-2$ and $K$ is an arbitrary field. Then the coordinate ring of the line $L$ is $S := K[t, u]$. 

The homogeneous ideal of a rope has been characterized in \cite{NNS}, Theorem 2.4. 
 
\begin{thm} \label{charact} Let $C \subset \pp^n$ be a curve of
degree at most $n-1$.
Then the following conditions are equivalent:
\begin{itemize}
\item[1.] $ C $ is an $(n-k)$--rope supported on the line $ L;$
\item[2.] $ I_C = ((I_L)^2, [x_0, \dots, x_r] B) $ where the matrix
$ B $ gives a map $ \varphi_B : \oplus_{j=1}^k S(-\beta_j-1) \to
S^{r+1}(-1) $ with $ \codim(I_k(B)) = 2;$
\item[3.] $ I_C = ((I_L)^2, F_1, \dots, F_k) $ where
$ V(F_1, \dots, F_k) \subset \pp^n$ is a scheme of codimension $
k$ which contains $L$ and is smooth at the points of $ L.$
\end{itemize}
\end{thm}

Actually, there is not only the matrix $B$ that plays an important role. 

\begin{obs} \label{rem-ass-matrices} {\rm
According to \cite{NNS}, Remark 2.7, every $(n-k)$--rope supported
on the  line $L$ is related to matrices $A, B$ with entries in $S$
such that there is an  exact sequence
\begin{equation} \label{vectfield}
0 \to \oplus_{j=1}^k \oo_{\pp^1}(-\beta_j-1)
\stackrel{B}{\lra} \oo_{\pp^1}^{r+1}(-1)
\stackrel{A}{\lra} \oplus_{i=0}^{r-k} \oo_{\pp^1}(\al_i-1)
\to 0
\end{equation}
where the ideals of maximal minors $ I_k(B) $ and $
I_{r+1-k}(A)$ have codimension 2. Note that $A^t$ is just the
syzygy matrix of $B^t$.}
\end{obs} 

For geometric characterizations of ropes we refer to \cite{NNS}, Remark 2.10. Now, we state a further description. 

\begin{lemma} 
A curve $C \subset \pp^n$ is a rope if and only if there is a line $L$ such that for every hyperplane $H \subset \pp^n$ containing $L$ there is an exact sequence 
$$
0 \to \cI_L (-1) \to \cI_C \to \cI_{C \cap H, H} \to 0. 
$$
\end{lemma} 

\begin{proof} 
If $C$ is a rope the condition is clearly satisfied. Conversely, the exact sequences provide $I_L^2 \subset I_C$. 
\end{proof} 

By the genus of a curve we always means its arithmetic genus. The genus of a rope is easily computed because in  \cite{NNS}, Lemma 2.6, Proposition 2.7,  and Corollary 2.8, it is shown. 

\begin{prop} \label{prop-genus} 
Using the notation above, an $(n-k)$--rope $C \subset \pp^n$ has genus 
$$
g = - \sum_{j=1}^k \beta_j = - \sum_{i=0}^{r-k} \alpha_i \leq 0. 
$$
Moreover, if $C$ is non-degenerate then $g \leq -k = \deg C - n$. 
\end{prop}

%%%%%%%%%%%%%%%%%%%%%%%%%%%%%%%%%%%%%%%%%%%%%%%%%%%%%%%%%%%%%%%%%%
\section{Minimal free resolutions of ropes on a line} \label{sec-res}

In this section, we determine a minimal free resolution of a rope
$ C $ supported on a line $ L$. Furthermore, we compute the
Hilbert function as well as the Hartshorne-Rao module $ M(C) $
of $ C.$
We use the notation of the previous section. Recall in particular
that $r = n-2$.

We begin by showing an interesting property of the generators of
the homogeneous ideal $ I_C $ of a rope $ C.$

\begin{lemma} \label{lem-ind}
If $ I_C = ((I_L)^2, F_1, \dots, F_k) $ is the
homogeneous  ideal of
an $ (n-k)$--rope $ C \subset \pp^n$ supported on the line $ L,$
then we have $$ ((I_L)^2, F_1, \dots, F_{h-1}) \cap (F_h) = (F_h)
\cdot I_L $$ for each $ h = 1, \dots, k.$
\end{lemma}

\begin{proof} Recall that $I_L = (x_0,\ldots,x_r)$.
According to Theorem \ref{charact} it is trivial that $ x_i F_h
\in (I_L)^2$ for every $ i=0, \dots, r,$ and so $ (F_h) \cdot I_L
\subseteq ((I_L)^2, F_1, \dots, F_{h-1}).$
Vice versa, assume $ G F_h \in ((I_L)^2, F_1, \dots, F_{h-1})$. We
may write $G$ as $ G = G' + G'' $ where $ G' \in I_L $ and $
G'' \in k[t,u]$. Since $ G' F_h \in (I_L)^2 $ we get $ G'' F_h \in
((I_L)^2, F_1, \dots, F_{h-1}) $. But the fact that $ G'' \in
k[t,u]$ implies $ G'' F_h \in (F_1, \dots, F_{h-1})$.  Since $
(F_1, \dots, F_h) $ is a regular sequence (cf.\ Theorem
\ref{charact}), we get that $ G'' = 0$. The claim follows.
\end{proof}

The previous proposition ensures that the sequences

\begin{equation*}
0 \lra I_L(-\beta_h-1) \stackrel{}{\lra} I_{h-1} \oplus (F_h)
\lra I_h \lra 0
\end{equation*}
are exact, where $ I_h = ((I_L)^2, F_1, \dots, F_h) $, $ h=1,
\dots, k,$ and $ I_0 = (I_L)^2$. Hence,  one can use the mapping
cone procedure to compute a free resolution of $ I_h $ from one of
$ I_{h-1}$ if one is able to find the comparison maps. Indeed, we
will carry out this program.
The starting point for the inductive procedure is a free
resolution of $ (I_L)^2$. The minimal free resolution of $
(I_L)^2$ is given by an Eagon--Northcott complex. However, we were
not able to determine the comparison maps between this
Eagon--Northcott complex and the Koszul complex which resolves $
I_L $. In order to overcome this problem we use a non-minimal free
resolution of $ (I_L)^2$. Probably, this resolution is known to specialists,
but we have no reference for it. In order to state it we need some
notation.
We write the Koszul complex $ P_{\bullet}$ which resolves  $ I_L $
as
\begin{equation*}\label{koszul}
P_{\bullet} \hspace*{2cm} 0 \lra \wedge^{n-1} P
\stackrel{\delta_{n-1}}{\lra} \dots \wedge^2
P \stackrel{\delta_2}{\lra} P \stackrel{\delta_1}{\lra} I_L \lra 0
\end{equation*}
where $ P = R^{n-1}(-1) = \oplus_{i=0}^r R e_i$ and $
\delta_1(e_i) = x_i $ for every i.

\begin{prop} \label{i2res}
The ideal $(I_L)^2 $ has the following
non-minimal free resolution
$$
\begin{array}{ccccccc}
& & \wedge^{n-1} P \otimes P & & & \wedge^2 P
\otimes P \\
0 & \lra & \oplus & \stackrel{\partial'_{n-1}}{\lra} & \dots &
\oplus & \stackrel{\partial'_2}{\lra} P \otimes P
\stackrel{\partial'_1}\lra (I_L)^2 \lra 0 \\[1ex]
 & & \wedge^{n-1} P & & & \wedge^2 P \\
 \end{array}
$$
where
$$
\partial'_i = \left (
\begin{array}{cc}
\delta_i \otimes \id_P & (-1)^i \partial_i \\ 0 & \delta_{i}
\end{array} \right ), \quad i = 3,\ldots,n-1,
$$
$$
\partial'_2 = (\delta_2 \otimes \id_P \quad \partial_2),
\hspace*{1cm}
\partial'_1 = \delta_1 \otimes \delta_1
$$
and $ \partial_i : \wedge^i P \to \wedge^{i-1} P \otimes P $ is
the canonical map defined by $$
\partial_i (u_1 \wedge \dots \wedge u_i) = \sum_{j=1}^i (-1)^{j+1} u_1 \wedge
\dots \wedge \hat{u_j} \wedge \dots \wedge u_i \otimes u_j. $$
\end{prop}

\begin{proof} Consider the following surjective map
$$ \al : I_L \otimes P \lra (I_L)^2 $$ defined by $ \al (x_i
\otimes  e_j) = x_i x_j$ for all $ i,j.$
The kernel of the map $ \al $ is generated by the trivial
relations $ x_i \otimes  e_j - x_j \otimes 
e_i.$ Moreover, if we embed $ I_L \otimes P $ in $ R \otimes P
\isom P $ we see that $ \ker \alpha $ is nothing else than $ \im
\delta_2 $. Therefore, we have the following exact sequence $$ 0
\lra \im \delta_2 \lra I_L \otimes P \lra (I_L)^2 \lra 0. $$
A free resolution of $ \im \delta_2 $ is provided by the Koszul
complex $P_{\bullet}$, while $ I_L \otimes P$ is resolved by
$P_{\bullet} \otimes P$.  It is easy to check that the following
diagram is commutative
$$
\begin{array}{ccccccccc} 0 & \lra & \im
\delta_2 & \lra & I_L \otimes P & \lra & (I_L)^2 & \lra & 0 \\
 & & \uparrow  & & \uparrow \\
 & & \wedge^2 P & \stackrel{\partial_2}{\lra} & P \otimes P \\
 & & \uparrow  & & \uparrow \\
 & & \vdots & & \vdots \\
 & & \uparrow  & & \uparrow \\
 & & \wedge^{n-1} P & \stackrel{\partial_{n-1}}{\lra} & \wedge^{n-2} P
\otimes P
\\
 & & \uparrow & & \uparrow \\
 & & 0 & & \wedge^{n-1} P \otimes P \\
 & & & & \uparrow \\
 & & & & 0 \end{array} \leqno(+)
 $$
Now the mapping cone procedure (cf., e.g., \cite{Weibel}) implies
our claim.
\end{proof}

In order to get a minimal free resolution of $ (I_L)^2 $ we have
to cancel redundant terms in the resolution above. For
stating this formally we need some more notation.

\begin{notation}
Since the map $ \partial_i : \wedge^i P \to \wedge^{i-1} P \otimes P
$ is split-injective we can identify $(\wedge^{i-1} P \otimes
P)/\partial_i(\wedge^i P)$ with the free $R$-module, say $D_{i-1}$, being a
direct summand of $\wedge^{i-1} P \otimes P$. Furthermore, the
previous diagram ($+$) shows that $\delta_{i-1} \otimes \id_P:
\wedge^{i-1} P \otimes P \to \wedge^{i-2} P \otimes P$ induces a
homomorphism
$$
d_{i-1}: D_{i-1} \to D_{i-2}, \quad i \geq 3.
$$
Similarly, putting $D_0 = R$ the map $\delta_{1} \otimes \id_P$
induces a homomorphism $d_1: D_1 \to (I_L)^2$.
Moreover, denote by
$$
\tau_B: \; Q := \oplus_{j=1}^k R(- \beta_j -1) \to P
$$
the extension of $\ffi_B$. Finally, define for $i \geq 2$
$$
\mu_{i}: \wedge^{i-1} P \otimes Q \to D_{i-1}
$$
to be the composition $\wedge^{i-1} P \otimes Q \to \wedge^{i-1} P \otimes P
\to D_{i-1}$ and let $\mu_1$ be the composition
$Q \stackrel{\ffi_B}{\lra} P \stackrel{\delta_1}{\lra} R$. 
\end{notation}
\medskip

The diagram ($+$) implies that
\begin{equation*}
0 \to D_{n-1}  \stackrel{d_{n-1}}{\lra} \dots \lra D_2 \stackrel{d_2}{\lra} D_1
\stackrel{d_1}{\lra} (I_L)^2 \to 0
\end{equation*}
is an exact sequence. Comparison with the Eagon--Northcott
complex shows that it is a minimal free resolution of $ (I_L)^2$.
This is the starting point for the description of the minimal free
resolution of an arbitrary rope on $ L.$

\begin{thm} \label{Cres}
Let $ C \subset \pp^n$ be an $ (n-k)$--rope
supported on the line
$ L $ with homogeneous ideal $ I_C = ((I_L)^2, [x_0,\dots,x_r] B).$
Then we have the following
 free resolution of $ I_C$:
\begin{equation*}
G_{\bullet} \hspace*{2cm} 0 \to G_n \stackrel{d_n'}{\lra} \dots \lra G_2 \stackrel{d_2'}{\lra} G_1
\stackrel{d_1'}{\lra} I_C \to 0
\end{equation*}
where
$$ G_i := D_i \oplus (\wedge^{i-1} P \otimes Q), 
% \cong
%(\wedge^i P \otimes P / \wedge^{i+1} P) \oplus
%(\oplus_{j=1}^k \wedge^{i-1} P(-\beta_j-1))
$$ 
and 
$$
d_i' := \left ( \begin{array}{cc}
d_i & (-1)^i \mu_i \\[1ex]
0 & \delta_{i-1} \otimes \id_Q 
\end{array} \right )  \qquad \mif 2 \leq i \leq n-1, 
$$ 
and 
$$
D_n := 0, \qquad d_1' := \left ( \begin{array}{lr}
d_1 & -\mu_1 \
\end{array} \right ), \quad d_n' := \left ( \begin{array}{c}
(-1)^n \mu_n \\[1ex]
\delta_{n-1} \otimes \id_Q 
\end{array} \right ).
$$
This resolution is minimal if and only if the rope $C \subset \pp^n$ is
non-degenerate.
\end{thm}

\begin{proof} We prove the claim by induction on $ k.$ For $ k=0$
the claim is shown by the discussion above. Notice that $G_n = 0$ in this
case.

Let $k \geq 1$. Write $I_C = ((I_L)^2, F_1,\ldots,F_k)$ and define
$\tilde{C} \subset \pp^n$ to be the $(n-k+1)$--rope defined by
$\tilde{I} = ((I_L)^2, F_1,\ldots,F_{k-1})  =
 ((I_L)^2, [x_0,\dots,x_r] \tilde{B})$, i.e.\ the matrix $\tilde{B}$
 is obtained from $B$ by deleting its last column. Applying
 the induction hypothesis
 to $\tilde{C}$ we denote
the maps and modules in the free resolution $\tilde{G_{\bullet}}$ of $\tilde{I}$
depending on $k$ by $\tilde{Q}, \tilde{\mu_i}$. Note that $Q =
\tilde{Q} \oplus R(- \beta_k -1)$. Denote by $\tau:  R(- \beta_k
-1) \to P$ the map such that $\tau_B = (\tau_{\tilde{B}},
\tau)$.

Using Lemma \ref{lem-ind} we have the short exact sequence
%(\ref{induction})
$$
0 \to I_L(-\beta_k -1) \to \tilde{I} \oplus (F_k) \to I_C \to 0.
$$
Define homomorphisms
$\gamma_i: \wedge^i P \to D_i \oplus (\wedge^{i-1} P \otimes \tilde{Q}),
\; i = 2,\ldots,n-1$ and
$\gamma_1: P \to D_1 \oplus  \tilde{Q} \oplus R(- \beta_k - 1)$ as
the sum of the composition
$$
\wedge^i P (- \beta_k - 1) \stackrel{\id_{\wedge^i P} \otimes \tau}{\lra}
\wedge^i P \otimes P \twoheadrightarrow D_i, \quad i = 1,\ldots,n-1,
$$
and the zero maps
$\wedge^i P (- \beta_k - 1) \to \wedge^{i-1} P \otimes \tilde{Q}$
and (for $i = 1$) the map $P (- \beta_k - 1)
\stackrel{\delta_1}{\longrightarrow} R(-\beta_k 
-1)$. A routine computation shows that $\gamma_{\bullet}$ is a
homomorphism between $P_{\bullet}$ and the direct sum of
$\tilde{G_{\bullet}}$ and the minimal free resolution of
$(F_k)$.
Hence the mapping cone procedure implies that $G_{\bullet}$ is a
free resolution of $I_C$. The description of the maps shows that
$G_{\bullet}$ is a minimal free resolution if and only if the
matrix $B$ does not contain non-zero elements of $K$ as entries,
which is equivalent to $C$ being non-degenerate (cf.\ Theorem
\ref{charact}).
\end{proof}

\begin{obs} {\rm
Let $C \subset \pp^n$ be a degenerate rope. Denote by $L \subset
\pp^n$  the smallest linear subspace containing $C$. Then
Theorem \ref{Cres} gives the minimal free resolution of $I_C$ as
module over the coordinate ring of $L$. From this resolution it is
easy to derive the minimal free resolution of $I_C$ as module over
the coordinate ring $R$ of $\pp^n$. In this sense Theorem
\ref{Cres} provides the minimal free resolution for {\it all} ropes
on a line.}
\end{obs}

From the free resolution $G_{\bullet}$ of $ I_C $ we easily derive some
consequences. We use the convention that $\binom{a}{b} := 0$ for
integers $a, b$ such that $a < b$.

\begin{corol} \label{hilbC}
Let $C \subset \pp^n$ be an $(n-k)$--rope C supported on a line.
Then we have:
\begin{itemize}
\item[(i)]
The Hilbert function $ h_C$ of $C$ is
$$
h_C(j) = \binom{j+1}{1}  + (r+1) \binom{j}{1}  - \sum_{h=1}^k \binom{j-\beta_h}{1}.
$$
\item[(ii)] The Castelnuovo-Mumford regularity of $C$ is $\max
\{\beta_j + 1 \s 1 \leq j \leq k \}$, i.e.\ it equals the maximal
degree of a minimal generator of $I_C$.
\end{itemize}
\end{corol}

\begin{proof} (i) follows from two easy computations:

(1) the Hilbert
function of the scheme associated to $ (I_L)^2 $ is the
alternating sum of the ranks of the modules of the Eagon-Northcott
complex, thus it is equal to $ \binom{j+1}{1} + (r+1) \binom{j}{1}$;

(2) the alternating sum of the ranks of  the
modules of the twisted Koszul complex $P_{\bullet} (-\beta_h)$ is equal to
$\binom{j-\beta_h}{1}$.

(ii) is a consequence of Theorem \ref{Cres} because the
Castelnuovo-Mumford regularity can be read off from a mimimal free
resolution.
\end{proof}

In \cite{NNS}, Proposition 3.1 we have computed the \HR\  of
a rope:

\begin{prop}\label{H1C} Let $ C \subset \pp^n$ be an $(n-k)$--rope
supported on the line
$ L $ with homogeneous ideal $ I_C=((I_L)^2, [x_0,\ldots,x_r] B).$ Then its
Hartshorne-Rao module $ H^1_*({\mathcal I}_C) $ is  as $R$--module
isomorphic to $ \coker(S^{r+1}(- 1) \stackrel{\varphi_A}{\lra}
\oplus_{j=0}^{r-k} S(\alpha_j -1))$.
\end{prop}

It is amusing to derive this result from Theorem \ref{Cres}.
Indeed, since $P$ has rank $n-1$ the last map $d_n'$ in the
resolution $\tilde{G}$ can be written as
$$
d_n' = \left ( \begin{array}{c}
(-1)^n \tau_B (-n+1) \\[1ex]
\delta_{n-1} \otimes \id_Q
\end{array}
\right ).
$$
It follows that
$$
H^1_*({\mathcal I}_C) \cong \ER^{n-1}(I_C, R)^{\vee} (n+1) \cong
(\coker (d_n')^*)^{\vee} (n+1) \cong (\coker \ffi_B^*)^{\vee} (2) \cong
\coker \ffi_A
$$
as claimed.
\smallskip

For later reference we state an immediate consequence of the
description of the Hartshorne-Rao module.

\begin{corol}\label{raoC} The Rao function
$ \rho_C : \ZZ \to \ZZ,\; \rho_C(i) = h^1_*({\mathcal I}_C (i)) $
of an $(n-k)$--rope
$ C \subset \pp^n $ supported on $ L $ is
$$
\rho_C(i) = \sum_{j=0}^{r-k} \binom{i + \alpha_j}{1} +
 \sum_{j=1}^k \binom{i-\beta_j}{1}  - (r+1) \binom{i}{1}.
$$
\end{corol}

%%%%%%%%%%%%%%%%%%%%%%%%%%%%%%%%%%%%%%%%%%%%%%%%%%%%%%%%%%%%%%%%%%%%%%%%%%%%%%%%%

\section{Global sections of the structure sheaf of a rope} \label{sec-struc-sheaf} 

In this section we compute the global sections of (all twists of) the structure
sheaf of a rope $ C $ supported on a line, we
describe the structure of $ H_*^0(C,\oo_C) $ both as $ S$--module
and as $ R$--module, and finally we determine the minimal free
resolution of $ H^0_*(C,\oo_C)$ as $ R$--module.

We continue using the notation of the previous sections.

\begin{obs} \label{rem-notation-gl-sec} {\rm
Given an $(n-k)$--rope $C \subset \pp^n$ we want to construct
$r+1-k = n-1-k$ global sections of its structure sheaf $\oo_C$
which are naturally associated to the matrix $A$ related to $C$
(cf.\ Remark \ref{rem-ass-matrices}).

We will construct these sections locally. To this end let $I =
\{i_1,\ldots,i_k\} \subset \{0,\ldots,r\}$ be a subset of $k$ distinct 
elements and put $J := \{0,\ldots,r\} \setminus I =
\{j_0,\ldots,j_{r-k}\}$ where we assume $j_0 < \ldots < j_{r-k}$.
We denote by $A_I$ the submatrix of $A$
formed by the columns of $A$ with index in $I$. Similarly, let
$B_I$ be the submatrix of $B$ consisting of the rows of $B$ with
index in $I$. The matrices $A_J, B_J$ are analogously defined.

Consider the open subset
$$
U_J := \{P \in \pp^n \s \det A_J(P) \neq 0\}.
$$
Since the determinants of $A_J$ and $B_I$ agree up to a sign (cf.\
\cite{BE-gor}) the matrix $A_J$ is invertible over $\oo_C(U)$ if and
only if $B_I$ has this property. 

Now we define local sections $z_0^J,\ldots,z_{r-k}^J \in \oo_C(U_J)$ by
$$
(z_0^J,\ldots,z_{r-k}^J) = (x_{j_0},\ldots,x_{j_{r-k}}) A_J^{-1}.
$$
Here we abuse notation by writing simply $\oo_C(U_J)$ instead of
$(\oo_C \otimes \oo_{\pp^n}(d))(U_J)$ for a suitable integer $d$.
Observe that $\deg z_i = 1 - \al_i \; (0 \leq i \leq r-k)$. The
next result shows that these local sections glue to global
sections.
}
\end{obs}

\begin{lemma} \label{lem-glueing} For all $i \in
\{0,\ldots,r-k\}$ and all subsets $J, J' \subset \{0,\ldots,r\}$
of $r+1-k$ elements
we have $z_i^J = z_i^{J'}$ in $\oo_C(U_J \cap U_{J'})$.
\end{lemma}

\begin{proof}
Since
$$
0 = A \cdot B = (A_I \ | \ A_J) \cdot \left (\begin{array}{cc}
B_I \\
B_J
\end{array} \right )
$$
% and $\det A_J = \pm \det B_I$ (cf.\ [Buchsbaum-Eisenbud paper???])
we get
$$
A_J^{-1} A_I = - B_J B_I^{-1}. \leqno(*)
$$
Hence $[x_0,\ldots,x_r] B \subset I_C$ implies
$$
x_I B_I + x_J B_J = [x_0,\ldots,x_r] B = 0
$$
in $\oo_C(U_J \cap U_{J'})$ where $x_I, x_J$ are the suitable rows
of variables. It follows
$$
\begin{array}{rcl}
x_I & = & - x_J B_J B_I^{-1} \\
& = & x_J A_J^{-1} A_I \quad ({\rm by} \ (*)) \\
& = & (z_0^J,\ldots,z_{r-k}^J) A_I \quad ({\rm by \ definition \  of} \ z_i^J).
\end{array}
$$
Putting $z^J = (z_0^J,\ldots,z_{r-k}^J)$ we have just shown $x_I =
z^J A_I$. The definition of $z^J$ implies $x_J = z^J A_J$. Hence we
obtain $[x_0,\ldots,x_n] = z^J A$ and in particular $x_{J'} = z^J
A_{J'}$. Plugging this into the definition of $z^{J'}$ we get
$$
z^{J'} = x_{J'} A_{J'}^{-1} = z^J A_{J'} A_{J'}^{-1} = z^J
$$
as claimed.
\end{proof}

\begin{notation} \label{notation-gl-sec}
Since the matrix $A$ has maximal rank at all points of the line $L
\cong \pp^1$ the open subsets $U_J$ cover the rope $C$ if $J$
varies in $\{0,\ldots,r-k\}$. Therefore, the local sections
$z_0^J,\ldots,z_{r-k}^J \in \oo_C(U_J)$ glue together to global
sections of $\oo_C$ which we denote by $z_0,\ldots,z_{r-k}$. Note
that $\deg z_i = 1 - \al_i$ (cf.\ the previous remark).
\end{notation}

\begin{obs} \label{rem-prop-gl-sec} {\rm
The end of the proof of the previous lemma shows
 that  we have in $H_*^0(C, \oo_{C})$
$$
[x_0,\ldots,x_r] = [z_0,\ldots,z_{r-k}] A.
$$
We will use this fact soon. }
\end{obs}

Using the global sections just defined we can describe all global
sections of $\oo_C(d)$ for all  $d \in \ZZ$. This
generalizes Migliore's result for double lines in $\pp^3$ (cf.\
\cite{mig}, Proposition 3.1).

\begin{thm} \label{thm-desc-gl-sec}
Let $ C \subset \pp^n$ be an $ (n-k)$--rope supported on the line $ L = \Proj
S$.
 Then its global sections of degree $d$ are
$$ H^0(C,\oo_C(d)) = \{ P + \sum_{i=0}^{r-k} Q_i z_i
\s P, Q_i \in S, \deg(P) = d, \deg(Q_i)=d+\al_i-1 \}.
$$
Furthermore, $H^0_*(C, \oo_C)$ is as $S$--module a free module
with basis $\{1, z_0,\ldots,z_{r-k}\}$.
\end{thm}

\begin{proof}
The first step is to  show that $1, z_0,\ldots,z_{r-k}$ generate a free
$S$--submodule of $H^0_*(C, \oo_C)$. Assume that there are
$f_{-1},f_0,\ldots,f_{r-k} \in S$ which are not all trivial such
that
$$
0 = f_{-1} + \sum_{i=0}^{r-k} f_i z_i
$$
in $H^0(C,\oo_C(d))$ for some $d$. We choose $J \subset
\{0,\ldots,r-k\}$ such that $U_J \neq \emptyset$. Then we have in
$\oo_C(U_J)$ using Remark
\ref{rem-prop-gl-sec}
$$
0 = f_{-1} + x_J \cdot A_J^{-1} \cdot \left ( \begin{array}{c}
f_0 \\
\vdots \\
f_{r-k}
\end{array} \right ).
$$
The exact sequence
$$
0 \to \ii_C(U_J) \to \oo_{\pp^n}(U_J) \to \oo_C(U_J) \to 0
$$
and the description of the homogeneous ideal of $C$ (Theorem
\ref{charact})
imply that there are $g_1,\ldots,g_k \in R_{(\det A_J)}$ such
that we have in $\oo_{\pp^n}(U_J)$
$$
\begin{array}{rcl}
f_{-1} + x_J \cdot A_J^{-1} \cdot \left ( \begin{array}{c}
f_0 \\
\vdots \\
f_{r-k}
\end{array} \right ) &  = & [x_0,\ldots,x_n] \cdot B \cdot
\left ( \begin{array}{c}
g_1 \\
\vdots \\
g_{k}
\end{array} \right ) \\
& = & x_I \cdot B_I \cdot \left ( \begin{array}{c}
g_1 \\
\vdots \\
g_{k}
\end{array} \right ) + x_J \cdot B_J \cdot \left ( \begin{array}{c}
g_1 \\
\vdots \\
g_{k}
\end{array} \right ).
\end{array}
$$
It follows $f_{-1} = 0$,
$$
B_I \cdot \left ( \begin{array}{c}
g_1 \\
\vdots \\
g_{k}
\end{array} \right ) = 0 \quad {\rm and} \quad A_J^{-1} \cdot \left
( \begin{array}{c}
f_0 \\
\vdots \\
f_{r-k}
\end{array} \right ) =  B_J \cdot \left ( \begin{array}{c}
g_1 \\
\vdots \\
g_{k}
\end{array} \right ).
$$
Since the matrix $B_I$ is invertible in $\oo_C(U_J)$ we get $f_0 =
\ldots = f_{r-k} = 0$, a contradiction. Thus, we have shown the second
claim.

In the second step we compare dimensions. Put $ V_d = \{ P + \sum Q_i z_i \s \deg(P) = d, \deg(Q_i) =  d+\al_i-1 \}$. Then we have
$$
\dim V_d = \binom{d+1}{1}  + \sum_{i=0}^{r-k} \binom{d+\al_i}{1}.
$$
But using Corollaries \ref{hilbC} and \ref{raoC} we get
$$
h^0(C, \oo_C(d)) = h_C(d) + \rho_C(d) =
\binom{d+1}{1}  + \sum_{i=0}^{r-k} \binom{d+\al_i}{1},
$$
i.e.\ $\dim V_d = h^0(C, \oo_C(d))$ for all $d$ and our claims
follow.
\end{proof} 

Later on we will deal with the $R$--module structure of $ H^0_*(C,
\oo_C)$. This structure is made partially explicit in the following statement. 

\begin{lemma} \label{lem-R-mod-str} Using the notation of Theorem
  \ref{thm-desc-gl-sec} we have for all $j = 0, 1, \ldots, r$
$$
x_j (P + \sum_{i=0}^{r-k} Q_i z_i) = \sum_{i=0}^{r-k} (P a_{ij}) z_i.
$$ 
\end{lemma} 

\begin{proof}
By construction of the global sections $z_l$  we have that $[x_0,\ldots,x_r] = [z_0,\ldots,z_{r-k}] A$. These relations specify the embedding $ R/I_C \lra H^0_*(C, \oo_C).$
But in $ R/I_C$ we have $ x_i x_j = 0$ for every $ i,j \leq r$. Thus, we obtain 
$$
0 = x_j \cdot [x_0,\ldots,x_r] = x_j \cdot[z_0,\ldots,z_{r-k}] A. 
$$
But the matrix $ A^t $
defines an injective map.  Hence, we get  $x_j \cdot [z_0,\ldots,z_{r-k}] = 0$ i.e. $ x_j z_l = 0 $ for every $ j =0, \dots, r $ and for
every $ l =0, \dots, r-k.$ Our claim follows. 
\end{proof} 

We close this section with the computation of the minimal
free resolution of $ H^0_*(C, \oo_C) $ as $ R$--module. In particular, we determine the maps. This also  completes the  description of the $R$--module structure of $ H^0_*(C, \oo_C)$. 

\begin{prop} \label{prop-res-struc} 
A graded minimal free  resolution of the $ R$--module $ H^0_*(C,
\oo_C) $ is 
\begin{eqnarray*}
 0 \longrightarrow \oplus_{i=-1}^{r-k} \wedge^{r+1} P(\al_i-1)
\stackrel{\varepsilon_{n-1}}{\longrightarrow} \dots \stackrel{\varepsilon_{2}}{\longrightarrow} \oplus_{i=- 1}^{r-k} P(\al_i-1) \stackrel{\varepsilon_{1}}{\longrightarrow} \hspace*{3.5cm}  \\
 \hfill  \oplus_{i=-1}^{r-k} R(\al_i-1) \stackrel{\varepsilon_{0}}{\longrightarrow} H^0_*(C, \oo_C) \longrightarrow 0 
\end{eqnarray*}
where $ \al_{-1}=1$ and, for $i \geq 1$,  
$$
\eeps_i = \left (\begin{array}{ccccc}
\delta_i & & & & \\
- \delta_i(A_{0, \bullet}) & \delta_i & & & \\
- \delta_i(A_{1, \bullet}) & & \delta_i & & \\
\vdots & & & \ddots & \\
- \delta_i(A_{r-k, \bullet}) & & & & \delta_i 
\end{array}
\right )
$$
Here, the non-specified parts are zero, $\delta_i$ denotes the $i$-th map in the Koszul complex with respect to $\{x_0,\ldots,x_r\}$ and $\delta_i(A_{j, \bullet})$ is the $i$-th map in the Koszul complex with respect to the elements of the $j$-th row of the matrix $A$. 
\end{prop}

\begin{proof} 
Let $ f_{-1}, f_0, \dots,f_{r-k} $ be the canonical basis of $ \oplus_{i=- 1}^{r-k}
R(\al_i-1)$. Then, the homomorphism 
$$ 
\eeps_0: \oplus_{i=-1}^{r-k} R(\al_i-1)
\lra H^0_*(C, \oo_C) 
$$ 
mapping $f_{-1}$ onto $1$ and $f_{i}$ onto the global section $z_i$ for $i \geq 0$ is clearly surjective. Using Remark \ref{rem-prop-gl-sec} and Lemma \ref{lem-R-mod-str} one checks that the sequence in the statement is a complex. The Buchsbaum-Eisenbud exactness criterion implies that (part of) it provides a minimal free resolution of $\im \eeps_1$. This allows to compute the Hilbert function of this module. Using Theorem \ref{thm-desc-gl-sec}, we see that the Hilbert functions of $\im \eeps_1$ and $\ker \eeps_0$ agree. It follows that $\im \eeps_1 = \ker \eeps_0$ completing the proof. 
\end{proof}

%%%%%%%%%%%%%%%%%%%%%%%%%%%%%%%%%%%%%%%%%%%%%%%%%%%%%%%%%%%%%%%%%%

\section{Global sections of the normal sheaf}
\label{sec-normal-sheaf}

In this section  we compute  global sections of  the
normal sheaf of a rope. We  obtain a complete description of these
global sections only for double lines. But the partial results for arbitrary ropes are powerful enough for our applications later on. 
The proofs use heavily the description
of the minimal free resolution of a rope $C$ as obtained in
Section \ref{sec-res}. Some results require particular care when the characteristic of the ground field is two.

It is well-known (cf., e.g., \cite{sernesi}) that the normal sheaf of $C$ is 
$$ H^0(C,{\mathcal N}_C) =
[\ker(H^0(C, \mbox{Hom}({\mathcal G}_1, \oo_C))
\stackrel{(d_2')^{*}}{\lra} H^0(C, \mbox{Hom}({\mathcal G}_2,
\oo_C))) ]_0
$$ 
where $ {\mathcal G}_1 $ and $ {\mathcal G}_2 $ are
the sheafifications of  $ G_1 $ and $ G_2,$ respectively, and the map $d_2'$ is defined in Theorem \ref{Cres}. 
Equivalently, the normal sheaf can be described as
\begin{equation}
H^0(C, {\mathcal N}_C) = [\ker(\mbox{Hom}(G_1, H^0_*(C, {\mathcal
O}_C)) \stackrel{(d_2')^{*}}{\lra} \mbox{Hom}(G_2, H^0_*(C,
{\mathcal O}_C)))]_0
\end{equation}
and we use this second formulation because it is more suitable for
the computation. 

We now fix some notation that we will use throughout this section. 
Let $C \subset \pp^n$ be an $(n-k)$--rope. 
Recall that $ G_1 = S_2(P) \oplus R(-\beta_1-1) \oplus \dots
\oplus R(- \beta_k-1)$ and that $P = R^{n-1}(-1) = \oplus_{i=0}^r R e_i$. Then we set for 
 $ \varphi \in \mbox{Hom}(G_1, H^0_*(C,
{\mathcal O}_C))$ (cf.\ Theorem \ref{thm-desc-gl-sec})
\begin{itemize}
\item $P^{ij} + \sum_{l=0}^{r-k} Q^{ij}_l z_l := \varphi(e_i e_j)  \qquad \mbox{where} \
P^{ij} \in [S]_2,\; Q^{ij}_l \in [S]_{\alpha_l+1} \;(0 \leq i, j \leq r);$
\item $P^s + \sum_{l=0}^{r-k} Q^s_l z_l := \varphi(f_s)  \qquad \mbox{where}  \ P^s \in [S]_{\beta_s+1},\; Q^s_l \in [S]_{\beta_s + \alpha_l} \; (1 \leq s \leq k)$, \\ 
and $ f_1, \dots, f_k $ is the canonical basis of $Q = R(-\beta_1-1)
\oplus \dots \oplus R(-\beta_k-1).$ (Notice that $P^s$ is just a polynomial with index $s$. It is not the $s$-th power of $P$. A similar comment applies to $Q_l^s$.) 
\end{itemize}

Furthermore, recall that $D_2 \cong \wedge^2 P \otimes P / \wedge^3 P$, \; 
$ G_2 = D_2 \oplus
P(-\beta_1-1) \oplus \dots \oplus P(-\beta_k-1) $ and that the map
$ d_2': G_2 \to G_1$ acts as follows
\begin{itemize}
\item $ d_2'([e_i \wedge e_j \otimes e_h]) = -x_i e_j e_h + x_j e_i e_h$ where we denote by $e_i e_j$ the class of $e_i \otimes e_j$ in $D_1$; 
\item $ d_2'(e^s_j) = \sum_{i=0}^r b_{is} e_i e_j + x_j f_s$ \qquad where $e_j^s := e_j \otimes f_s$ denotes the $(j+1)$-st basis element of $P(-\beta_s - 1)$. 
\end{itemize}
Then we obtain for the map $(d_2')^{*}(\varphi) = \ffi \circ d_2': G_2 \to H^0_*(C,
{\mathcal O}_C)$ using Lemma \ref{lem-R-mod-str}
\begin{itemize}
\item $ (d_2')^{*}(\varphi)([e_i \wedge e_j \otimes e_h]) =
\sum_{l=0}^{r-k} (- P^{jh} a_{li} + P^{ih} a_{lj}) z_l;$
\item $ (d_2')^{*}(\varphi)(e^s_j) = \sum_{i=0}^r b_{is} P^{ij} +
\sum_{l=0}^{r- k} (\sum_{i=0}^r b_{is} Q^{ij}_l + P^s a_{lj})
z_l.$
\end{itemize}

Combined with Theorem \ref{thm-desc-gl-sec} we get the following result. 

\begin{prop} \label{prop-lin-systems}
Adopting the notation above we have $ \varphi \in H^0(C, {\mathcal N}_C) $ if and only if the
following conditions are satisfied
\begin{itemize}
\item[1.] $ P^{ih} a_{lj} - P^{jh} a_{li} = 0 \qquad \mbox{for every } l=0,
\dots, r-k, \quad i,j,h = 0, \dots, r;$ \\[-8pt]
\item[2.] $ \sum_{i=0}^r b_{is} P^{ij} = 0 \qquad \mbox{for every } s=1, \dots,
k, \quad j=0, \dots, r;$ \\[-8pt]
\item[3.] $ \sum_{i=0}^r b_{is} Q^{ij}_l + P^s a_{lj} = 0 \qquad \mbox{for
every } s=1, \dots, k, \quad l = 0, \dots, r-k, \quad j = 0,
\dots, r.$
\end{itemize} \end{prop} 

\begin{obs} \label{rem-free}
\rm (i) The conditions above are linear systems over the ring $S$ with unknowns $ P^{ij}$'s (conditions 1. and 2.) and with unknowns $ P^s$'s and $ Q^{ij}_l$'s (condition 3.). 

(ii) In the conditions 1.\ - 3.\  the $ Q_l^s $'s appearing in the definition of $\ffi(f_s)$ are not involved. Thus, they are free parameters in the description of the global sections of ${\mathcal N}_C$. 
\end{obs}

Next, we will solve the linear systems given by 1.\ and 2.\ above. In case of degenerate ropes we need some further notation. 

\begin{notation} \label{not-deg-rope}
Let $C \subset \pp^n$ be a rope spanning a $\pp^m$. Then $C$ is degenerate if and only if $m < n$. In this case, we may and will assume that the matrices $A, B$ associated to $C$ have the shape 
$$
A = (A'\;\ 0)
$$
and 
$$
B = \left ( \begin{array}{cc}
B' & 0 \\
0 & I_{n-m} 
\end{array} \right )
$$
where $A'$ and $B'$ are the matrices associated to the non-degenerate rope $C$ considered as subscheme of $\pp^m = \mbox{Proj} (K[x_0,\ldots,x_{m-2}, t, u])$, thus $\pp^m = \{x_{m-1} = \ldots = x_{n-2} = 0\}$. 
\end{notation} 

Now, we are ready for solving systems 1.\ and 2. 

\begin{prop}\label{pij} If $ C \subset \pp^n$ is an $(n-k)$-rope on $ L,$ then we have for $\ffi \in H^0(C, {\mathcal N}_C)$ 
\begin{enumerate}
\item[(a)] $ (P^{00},P^{01},P^{11}) = c
  (a^2_{00},a_{00}a_{01},a^2_{01}) $ for some $ c \in K $ and $P^{i j}
  = 0$ for all other indices  in case  $ C$ has degree $ 2 $ and
  genus 
$ -1;$
\item[(b)] $ P^{ij} = 0 $ for every $ i,j = 0, \dots,
r,$ otherwise.
\end{enumerate}
\end{prop}

\begin{proof}  First, we assume that $C$ is non-degenerate. We distinguish two cases. 

\it Case 1: \rm Let $ \deg(C) > 2.$ 

Assume on the contrary that $P^{j h} \neq 0$ for some $j, h$. 
We denote by $ {\bf a}_0, \dots, {\bf a}_r $ the columns of the
matrix $ A$. Then part of the system 1. can be written as 
$$ 
P^{ih} {\bf
a}_j = P^{j h} {\bf a}_i \qquad \mbox{for every } i,h = 0, \dots,
r.
$$
Since $ P^{jh} \not= 0 $  we get that $ P^{ih} \not= 0
$ for every $ i$ because every column $ {\bf a}_i$ is not zero. The latter is true since otherwise one $\beta_k$ would vanish contradicting the assumption that 
$C $ is non-degenerate (cf.\ Theorem \ref{charact}). It follows that the matrix $ A $ has rank 1 at a general point of $ L,$ but rank $ A = r+1-k > 1.$ This contradiction shows that $P^{ij} = 0 $ for every $ i,j = 0, \dots, r.$ 

\it Case 2: \rm Let $ \deg(C) = 2.$

In this case we consider the linear system 2. It can be
written as 
$$ (P^{0j}, \dots, P^{rj}) \cdot B = 0 \qquad \mbox{for every
} j = 0, \dots, r.
$$ 
Since $A^t$ is the syzygy matrix of $B^t$ we get $ (P^{0j}, \dots, P^{rj}) = c_j A$ for some $c_j \in S$. Notice that in case $\deg C = 2$ we have $A = (a_{0 0},\ldots,a_{0 r})$. 
Moreover, recalling that $ P^{ij} = P^{ji} $ we see that $ c_i
a_{0j} = c_j a_{0i},$ for each $ i,j.$
Note that every entry of $A$ has degree $\alpha_0 \geq 1$. We distinguish three 
cases: 

$ (1) \qquad \alpha_0 \geq 3;$ 

$ (2) \qquad \alpha_0 = 2;$ 

$ (3) \qquad \alpha_0 = 1.$ 

We use the fact that every $ P^{ij}$ has degree 2. 
It implies in case (1) that every $ c_i $ is zero for degree reasons, thus  $P^{ij} = 0$ for every $ i,j.$ 

In case (2) the $c_j$'s are in $K$.  If one of them is not zero then all the entries of $ A $ differ by a scalar only.  Thus, $A$ has degree zero syzygies, i.e.\ $ C $ is
degenerate. This contradiction shows  the claim  in this case, too. 

In 
case (3) the rope $C$ has genus $-1$, thus it is degenerate if $n \geq 4$ (cf.\ \ref{prop-genus}). But if $n = 3$, the system 1.\ and  the system 2.\ read as 
$$ 
P^{00}a_{01} - P^{01}a_{00} = 0 \qquad \mbox{and} \qquad P^{01}a_{01} -
P^{11}a_{00} = 0 
$$ 
because $ A = (a_{00}, a_{01}) $ and $ B =
(a_{01}, -a_{00}) $ where $ a_{00} $ and $ a_{01} $ are two independent linear
forms.
Every solution of that system has the form 
$$ (P^{00},P^{01},P^{11}) = c
(a^2_{00}, a_{00}a_{01}, a^2_{01})  \qquad \mbox{for some } c \in
K  
$$ 
completing the proof for non-degenerate ropes.  

Second, we assume that $C$ is degenerate. Then using the system 2.\ and the block structure of its matrix $B$ (cf.\ \ref{not-deg-rope}), we can reduce to the case of non-degenerate ropes. 
\end{proof}

Having resolved the systems  1.\ and 2.\ it remains to consider the linear system given by condition 3.\ in Proposition \ref{prop-lin-systems}. We are not able to do this in full generality except in the case of double lines, i.e.\ $2$--ropes. We use Notation \ref{not-deg-rope}. 
 
\begin{prop} \label{prop-sys-3} 
Let $ C $ be a $ 2$--rope of genus $g$. Then the polynomials $P^s, Q_0^{i j}$ form a solution of the linear system 3.\ if and only if 
\begin{itemize} 
\item[(i)] in case $\chr K \neq 2$: \\ 
$$ 
P^s = \left \{\begin{array}{ll} 
-\frac{1}{2} \sum_{j=0}^{m-2} \lambda_j b_{js} & \mif  1 \leq s  \leq m-2 \\[1ex]
- \lambda_s & \mif m-1 \leq s \leq n-2 = r 
\end{array} \right. 
$$ 
and 
$$ 
Q^{ij}_0 = \left \{ \begin{array}{ll} 
\frac{1}{2} (\lambda_i a_{0j} + \lambda_j a_{0i}) & \mif  0 \leq i, j \leq  m-2 \\[1ex]
\lambda_j a_{0 i} & \mif 0 \leq i \leq m-2 < j \leq r \\[1ex]
0 & \mif m-1 \leq i, j \leq r
\end{array} \right.
$$ 
where $\lambda_0,\ldots,\lambda_r \in S$ are any linear forms; 
\item[(ii)] in case $\chr K = 2$:  \\
$$
Q_0^{i j} = \left \{ \begin{array}{ll} 
0 & \mif 0 \leq i = j \leq m-2 \\[1ex]
\sum_{h=1}^{m-2} \adj((B'_j)^t)_{i h} P^k & \mif 0 \leq i, j \leq m-2\; \mbox{and}\;i \neq j \\[1ex]
a_{0 i} P^j & \mif 0 \leq i \leq m-2 < j \leq r \\[1ex]
0 & \mif m-1 \leq i, j \leq r
\end{array} \right. \quad \mbox{provided}\; g \leq -2
$$ 
whereas 
$$
Q_0^{i j} = \left \{ \begin{array}{ll} 
c' a_{0 j}^2 & \mif i = j = 0, 1 \\[1ex]
c' a_{0 j} a_{0 i} + a_{0 1} P^1 & \mif i = 0, j = 1\\[1ex]
0 & \mbox{otherwise} 
\end{array} \right. \; \mbox{for some}\ c' \in K\; \mif g = -1.
$$
Here $\adj(\_)$ denotes an adjoint matrix, $B'_j$ is the matrix obtained from $B'$ by deleting its row with index $j$ and $M_{i h}$ is the entry at position $(i, h)$ of the matrix $M$.\footnote{For convenience, the subscripts of the rows of $ B $ range 
from $ 0 $ to $ r.$}
\end{itemize} 
\end{prop}

\begin{proof} At first, we state some notation.
$ C $ being  a $ 2$--rope we have $ k = r.$ Moreover, $ A $ is a
matrix of type $ 1 \times (r+1),$ while $ B $ is a matrix of type
$ (r+1) \times r,$ and they are related (using the Hilbert-Burch theorem) by 
$$ 
a_{0j} = (-1)^j \det(B_j). 
$$ 
Observe that all $a_{0 j}$ are not zero because $C$ is assumed to be non-degenerate. 

As last piece of notation, we set $ {\bf Q}^j =
(Q^{0j}_0, \dots, Q^{rj}_0)^t $ and $ {\bf P} = (P^1, \dots,
P^r)^t.$ Note that $ Q^{ij}_0 = Q^{ji}_0$ by their definition.
Hence, the linear system 3. can be written as
\begin{equation}\label{sys}
B^t {\bf Q}^j + a_{0j} {\bf P} = {\bf 0}, \qquad j=0, \dots, r.
\end{equation} 
Since $a_{0 j} = 0$ if $m-1 \leq j \leq n-2 = r$ (cf.\ Notation \ref{not-deg-rope}) we get for $ j = m-1,\ldots,r$ 
$$
\left ( \begin{array}{cc}
B' & 0 \\
0 & I_{n-m} 
\end{array} \right ) \cdot {\bf Q}^j = 0. 
$$ 
The last $n-m$ equations immediately imply $\left ( \begin{array}{c}
Q^{m-1, j}_0  \\
\vdots \\
Q^{r, j}_0 
\end{array} \right ) = 0$. 
In particular, symmetry $Q^{i, j}_0 = Q^{j, i}_0$ holds if $m-1 \leq i, j \leq r$. 

The first $m-1$ equations of the last system above read as 
$$
(B')^t \cdot \left ( \begin{array}{c}
Q^{0, j}_0  \\
\vdots \\
Q^{m-2, j}_0 
\end{array} \right ) = 0. 
$$
Since $(A')^t$ is the syzygy matrix of $(B')^t$ it follows 
$$
\left ( \begin{array}{c}
Q^{0, j}_0  \\
\vdots \\
Q^{m-2, j}_0 
\end{array} \right ) = \lambda_j (A')^t \qquad (m-1 \leq j \leq r) 
$$
for some $\lambda_j \in [S]_1$ by comparing degrees. 

For $j = 0,\ldots,m-2$ system (\ref{sys}) provides 
$$
(B')^t \cdot \left ( \begin{array}{c}
Q^{0, j}_0  \\
\vdots \\
Q^{m-2, j}_0 
\end{array} \right ) + a_{0 j} (B')^t \cdot \left ( \begin{array}{c}
P^1 \\
\vdots \\
P^{m-2}
\end{array} \right ) = 0 \; \mbox{and}\; \left ( \begin{array}{c}
Q^{m-1, j}_0  \\
\vdots \\
Q^{r, j}_0 
\end{array} \right ) = - a_{0 j} \left ( \begin{array}{c}
P^{m-1} \\
\vdots \\
P^{n-2}
\end{array} \right ). 
$$
Hence, symmetry implies if $m-1 \leq i \leq r$ and $0 \leq j \leq m-2$ 
$$
\lambda_i a_{0 j} = Q^{j i}_0 = Q^{i j}_0 = - a_{0 j} P^i. 
$$
Since by assumption $C$ spans $\pp^m$ we must have $a_{0 j} \neq 0$ if $0 \leq j \leq m-2$. It follows that 
$$
P^i = - \lambda_i \qquad \mbox{and} \qquad Q^{i j}_0 = -a_{0 i} P^j \qquad \mif m-1 \leq i \leq r\; \mbox{and}\; 0 \leq j \leq m-2. 
$$

It remains to show the claims for the indices $i, j \leq m-2$, i.e.\ to consider the non-degenerate rope in $\pp^m$ associated to the matrices $A', B'$. In order to simplify notation we assume that $C$ itself is non-degenerate, i.e.\ $n = m,\; A' = A,\; B' = B$.  

Now, we will show that system (\ref{sys})  can be solved if we also consider {\bf P} as given and only the ${\bf Q}^{j}$'s as unknowns. To this end we first determine a particular solution of the system (depending on {\bf P}) and then all the solutions of the associated  homogeneous system. 

$ \bullet $ Particular solution of (\ref{sys}). 

We will describe a solution satisfying $Q_0^{j j} = 0$. 
We  denote by $ {\bf \overline Q}^j $ the vector
obtained from $ {\bf Q}^j $ by deleting the entry $Q_0^{j j}$. Then, assuming $Q_0^{j j} = 0$, 
the system (\ref{sys}) becomes 
$$ B_j^t \cdot {\bf \overline Q}^j = - a_{0 j} P = 
-(-1)^j \det(B_j) {\bf P}.
$$ 
If we multiply both sides by $\mbox{adj}(B_j^t)$  we can
cancel $ \det(B_j^t) $ because it is not zero (thus it is not
a zero divisor). Hence,  we get $ {\bf \overline Q}^j = (-1)^{j+1}
\mbox{adj}(B_j^t) {\bf P}.$ Therefore, a particular solution $ {\bf \tilde Q}^j$ is given
by inserting a zero entry after the $ j$-th entry of $ {\bf \overline Q}^j$. 

$ \bullet $ Solution of the associated homogeneous system $B^t {\bf Q}^j = 0,\; j = 0,\ldots,r$. 

The syzygy argument above shows that every solution is of the form $ {\bf Q}^j = \lambda_j A^t,$
where $ \lambda_j \in [S]_1$. 
\smallskip 

Hence, the general solution of (\ref{sys}) is
\begin{equation}\label{solsys}
{\bf Q}^j = {\bf \tilde Q}^j + \lambda_j A^t, \qquad j = 0,\ldots,r. 
\end{equation}
But, as noted above,  the polynomials $Q_0^{i j}$ also satisfy $ Q^{ij}_0 = Q^{ji}_0$ for all $i, j \leq r$. For $i=0$ this extra requirement reads as 
$$
-[\mbox{adj}(B_0^t)]_{j-\mbox{row}} {\bf P} + (-1)^j \lambda_0
\det(B_j) = (-1)^{j+1} [\mbox{adj}(B_j^t)]_{0-\mbox{row}} {\bf
P} + \lambda_j \det(B_0).  
$$ 
Thus, we get 
$$ 
\{(-1)^j
[\mbox{adj}(B_j^t)]_{0-\mbox{row}} - [\mbox{adj}(B_0^t)]_{j-
\mbox{row}} \} {\bf P} = \lambda_j \det(B_0) - (-1)^j \lambda_0
\det(B_j).
$$
If we denote by $ B_{ij;h} $ the submatrix obtained from $ B $ by
deleting rows $ i $ and $ j $ and the column $ h$,  we obtain 
$$
\lambda_j \det(B_0) - (-1)^j \lambda_0 \det(B_j) = (-1)^j
\sum_{h=1}^r (- 1)^h \det(B_{0j;h}) (\lambda_j b_{jh} + \lambda_0
b_{0h}) 
$$ 
and 
\begin{eqnarray*}
(-1)^j [\mbox{adj}(B_j^t)]_{0-\mbox{row}} -
[\mbox{adj}(B_0^t)]_{j- \mbox{row}}  = \hspace*{4.5cm} \\
\hfill  2(-1)^j [\det(B_{0j;1}), -\det(B_{0j;2}), \dots, (-1)^{r-1}
\det(B_{0j;r})]. 
\end{eqnarray*}
 Substituting it follows  
\begin{eqnarray} \label{eqn-dist} \nonumber
 2[(-1)^{j+1}
\det(B_{0j;1}), \dots, (-1)^{j+r}\det(B_{0j;r})] \cdot {\bf P} = \hspace*{4.5cm} \\ 
\hfill \sum_{h=1}^r (-1)^{j+r+1} \det(B_{0j;h}) (\lambda_j b_{jh} +
\lambda_0 b_{0h}), \qquad j=1, \dots, r.
\end{eqnarray} 
Now, we distinguish according to the characteristic. 

{\it Case 1}: Assume $\chr K \neq 2$. 

Using matrix notation the previous equations  read as  
$$ 
2 \mbox{adj}(B_0^t) {\bf P} =
{\bf n} 
$$ 
where $ {\bf n} $ is defined in the obvious way. Multiplying
both sides by $ B_0^t $ provides 
$$ 
2 a_{00} {\bf P} = B_0^t {\bf
n}.
$$
If we compute the $s$-th entry of both sides we obtain  
$$ 2 a_{00}
P^s = \sum_{j=1}^r b_{js} \sum_{h=1}^r (-1)^{h+j+1}
\det(B_{0j;h})(\lambda_j b_{jh} + \lambda_0 b_{0h}) = -a_{00}
\sum_{j=0}^r \lambda_j b_{js}.
$$ 
Being $ \mbox{char}(K) \not= 2,$
it follows  that 
$$ P^s = -\frac 12 \sum_{j=0}^r \lambda_j b_{js},
\qquad s=1, \dots, r.
$$
If we plug these formulas into (\ref{solsys}), we get 
$$ Q_0^{ij} =
Q_0^{ji} = \frac 12 (\lambda_j a_{0i} + \lambda_i a_{0j}) 
$$ 
as claimed. 

{\it Case 2}: Assume $\chr K = 2$. 

Then (\ref{eqn-dist}) shows that the symmetry conditions $Q_0^{i j} = Q_0^{j i}$ are equivalent to 
$$
\lambda_j a_{0 i} = \lambda_i a_{0 j} \qquad \mbox{for all}\; i, j =
0,\ldots,r.   
$$
These equations are similar to the system 1.\ considered in
Proposition \ref{pij}. Analogous arguments provide   
$\lambda_i = 0$ for all $i$ if $g \leq -2$. If $g = -1$ we get $\lambda_i = c' a_{0 i}$ for $i = 0, 1$ and some $c' \in K$ and $\lambda_i = 0$ if $ 2 \leq i \leq r$. Using equations (\ref{solsys}) the asserted necessary conditions follow. 

In both cases it is easy to see that the stated conditions are sufficient as well.  
\end{proof}

Thanks to the previous computations, we obtain the following result. 

\begin{corol}\label{h0nc2} 
Let $ C \subset \pp^n\; (n \geq 3)$ be a double line of genus $g$.  Then we have  $$ 
h^0(C,{\mathcal N}_C) = \left\{
\begin{array}{cl} 
4 (n-1) & \mif\; g = -1   \\
(n-1)(3-g)-1 & \mif\; \chr K \neq 2 \ \mbox{and} \ g \leq -2 \\
n (3-g) - 6 & \mif\;  \chr K = 2 \ \mbox{and} \ g \leq -2
 \end{array} \right. 
$$
\end{corol}

\begin{proof} 
It is clear that the two last results give not only necessary, but also sufficient conditions for the global sections of ${\mathcal N}_C$. Thus, 
$ h^0(C,{\mathcal N}_C) $ equals the number of free
parameters in the description of the solutions of the systems 1.\ - 3.\ in Proposition \ref{prop-lin-systems}. 
Counting these parameters is easy. 

{\it Case 1}: Assume $\chr K \neq 2$. 

Then we have 
\begin{itemize}
\item $ r+1 $ arbitrary polynomials $ \lambda_0, \dots, \lambda_r $ in $ [S]_1$ (cf.\ Proposition \ref{prop-sys-3}) 
giving $ 2(r+1) $ free parameters;
\item $r$  arbitrary polynomials $ Q_0^1, \dots, Q_0^r $ in $ S $ of degrees $
\alpha_0 + \beta_1, \dots, \alpha_0 + \beta_r,$ respectively, (cf.\ Remark \ref{rem-free}) providing $ r(\alpha_0 + 1) + \alpha_0 $ free parameters because $\beta_1 + \ldots + \beta_r = \alpha_0 = -g$ (cf.\ Proposition \ref{prop-genus}). 
\end{itemize}
Moreover, if $ g \leq -2$ then $ P^{ij} = 0 $ for every $ i,j,$
while, if $ g = -1 $  we have one more free parameter (cf.\ 
Proposition \ref{pij}). The claim follows. 

{\it Case 2}: Assume $\chr K = 2$. 

Then we have in case $g \leq -2$ 
\begin{itemize} 
\item arbitrary polynomials $P^i \in [S]_{\beta_i + 1}$ (cf.\ Proposition \ref{prop-sys-3}) giving $2 r - g$ parameters; 
\item as in Case 1 arbitrary polynomials $ Q_0^1, \dots, Q_0^r $ in $ S $ of degrees $
\alpha_0 + \beta_1, \dots, \alpha_0 + \beta_r,$ respectively.  
\end{itemize} 
The claim follows. 

If $g = -1$ we have 
\begin{itemize}
\item an arbitrary constant $c' \in K$ and an arbitrary $P^1 \in [S]_2$ by Proposition \ref{prop-sys-3}; 
\item an arbitrary $c \in K$ due to Proposition \ref{pij}; 
\item an arbitrary $Q_0^1 \in [S]_2$. 
\end{itemize} 
Summing up, we get again $h^0(C, {\mathcal N}_C) = 4 (n-1)$. 
\end{proof}

Now, we want to exhibit a family of solutions of the linear system
3.\ for an arbitrary rope $ C $ of degree $n-k \geq 3$ that satisfy an additional condition. 

\begin{condition} \label{cond} 
Assume that the vector $ {\bf P} = (P^1,\ldots,P^k)$ is a combination of the rows of $ B,$ i.e.\ 
$$ 
{\bf P} = B^t \cdot {\bf \lambda} 
$$ where $ {\bf
\lambda} = \rm (\lambda_0, \dots, \lambda_r)^t$ and $\lambda_i \in S$. Note that $
\deg(\lambda_i) = 1 $ for each $ i$ by comparing degrees. 
\end{condition} 

We are going to determine all solutions of the system 3.\ satisfying this condition. Observe that the condition is automatically satisfied if $C$ is a $2$--rope and $\chr K \neq 2$ according to  Proposition \ref{prop-sys-3}.

For the remainder of this section we assume that the rope $C$ has degree at least 3. Setting $ {\bf Q}_l^j = (Q_l^{0,j}, \dots, Q_l^{rj})^t$ we can
write the system 3. as 
$$ 
B^t {\bf Q}_l^j + a_{lj} B^t {\bf
\lambda} = {\bf 0}, \qquad l=0, \dots, r-k, \; j = 0,\ldots,r, 
$$ 
thus 
$$
B^t ({\bf Q}_l^j + a_{lj} {\bf \lambda}) = {\bf 0}.
$$
Using Remark \ref{rem-ass-matrices} again, we get 
for each $ l, j$  
$$ 
{\bf Q}_l^j +
a_{lj} {\bf \lambda} = A^t {\bf c}^{lj} 
$$ where $ {\bf c}^{lj} =
(c_0^{lj}, \dots, c_{r-k}^{lj})^t $ and $ \deg(c_h^{lj}) =
\alpha_l - \alpha_h + 1$\ for $ h = 0,
\dots, r-k$ and $ \alpha_l -
\alpha_h + 1 < 0 $ means that $ c_h^{lj} = 0.$ 
Since $ Q_l^{ij} = Q_l^{ji}$,  the parameters $ c_h^{lj} $
satisfy the following relations 
$$ 
- \lambda_i a_{lj} + \sum_{h=0}^{r-k} a_{hi} c_h^{lj} = -
\lambda_j a_{li} + \sum_{h=0}^{r-k} a_{hj} c_h^{li}, \qquad  i,j = 0, \dots, r,   
$$ 
or equivalently
\begin{equation} \label{sisc} 
\sum_{h=0}^{r-k} a_{hi} c_h^{lj} -
\sum_{h=0}^{r- k} a_{hj} c_h^{li} = \lambda_i a_{lj} - \lambda_j
a_{li} \qquad \mif 0 \leq j < i \leq r\; \mbox{and}\; 0 \leq l \leq r-k.
\end{equation} 

We think of the  system (\ref{sisc}) as a linear
system with unknowns $ c_h^{lj}$.  It remains to compute its solutions. 

It is easy to see that particular solution of (\ref{sisc}) is 
$$ 
c_h^{lj} = \left\{ \begin{array}{cccc} 0 & \  & \mbox{if } & h
\not= l \\ - \lambda_j & \  & \mbox{if } & h = l \end{array}
\right. \qquad  (j = 0, \dots, r).
$$
Now, we consider the homogeneous linear system associated to
(\ref{sisc}), i.e.\ 
\begin{equation} \label{hsisc} \sum_{h=0}^{r-k} a_{hi} c_h^{lj} -
\sum_{h=0}^{r- k} a_{hj} c_h^{li} = 0 \qquad 0 \leq j < i \leq r.
\end{equation}
The coefficient matrix with respect to the variables $ c_h^{l0}, c_h^{l1},
\dots, c_h^{lr} $ is 
$$ 
M_h = \left[ \begin{array}{ccccc} a_{h1} &
-a_{h0} & 0 & \dots & 0 \\ a_{h2} & 0 & -a_{h0} & \dots & 0 \\
\vdots & \vdots & \  & \ddots & \  \\ a_{hr} & 0 & \  & \dots &
-a_{h0} \\ 0 & a_{h2} & -a_{h1} & \dots & 0 \\ \vdots \\ 0 & \dots
& 0 & a_{hr} & -a_{h,r-1} \end{array} \right].
$$
Notice that $ M_h $ is the matrix that describes the map $ F
\to \wedge^2 F $ of the Koszul complex $ K(a_{h0}, \dots, a_{hr})
$ associated to the  row of $A$ with index $h$. 

The whole coefficient matrix (\ref{hsisc}) is $ M = (M_0\, |\, M_1\, |\, \dots\, |\,
M_{r-k}).$ Thus, it is easy to verify that the following vectors
are solutions of (\ref{hsisc}): 
$$ (\dots, 0, A_h, 0, \dots)^t, 
\qquad h=0, \dots, r-k,  
$$ 
$$ 
(\dots, 0, A_{h'}, 0, \dots, 0, A_h, 0,
\dots)^t, \qquad 0 \leq h < h' \leq r-k 
$$ 
where $ A_0, \dots, A_{r-k}$ are the rows of $ A. $ 
These solutions are independent over $S$ because the rows of $A$ have this property. The latter is true since the maximal minors of $A$ define an ideal of codimension two. Hence, we have found $\binom{r+2-k}{2}$ independent solutions of (\ref{hsisc}).

Now, we want to estimate the rank of $ M$.
To this end we set $ M_h = ( m_{h0}, m_{h1}, \dots, m_{hr})$
where $ m_{hj} $ denotes the $j$-th column of $ M_h$. Consider
the following submatrix of $ M $ 
\begin{eqnarray*}
 M' := ( m_{00}, \dots,
m_{r-k,0}; \dots; m_{0k}, \dots, m_{r-k,k}; \hspace*{4cm}\\ 
\hfill m_{0,k+1},
\dots, m_{r-k-1,k+1}; \dots; m_{0,r-2}, m_{1,r-2}; m_{0,r-1})  
\end{eqnarray*}
that is obtained by picking the first $ r $ columns of $ M_0, $ the
first $ r-1 $ columns of $ M_1, \dots, $ and the first $ k $
columns of $ M_{r-k}$ and reordering them suitably. $ M'
$ is an upper triangular  block matrix. If we delete the first $ k-1
$ rows of the first block, the first $ k-2 $ rows of the
second block, $ \dots,$ the first row of the $ (k-1)$--st block,
we obtain a square block matrix with square blocks on the diagonal that is 
upper triangular, has order $ (r+k)(r+1-k)/2$, and  
determinant 
$$ 
a_{0r} \det(A_{0,1;r-1,r}) \cdots
\det(A_{0,\dots,r- 1-k;k+1, \dots, r}) [\det(A_{0, \dots, r-k; k,
\dots, r})]^k 
$$ 
where $ A_{i_1, \dots, i_p; j_1, \dots, j_p} $ denotes 
the submatrix of $ A $ obtained by picking the rows $ i_1, \dots,
i_p $ and the columns $ j_1, \dots, j_p.$ Analogously, there are other 
 square submatrices of $ M $ of the same order
such that their determinants are products of the $k$-th power of a maximal minor
of $ A$ and other subminors of it and all such possibilities occur.  But $A$ has maximal rank. It follows that  the rank of $ M $ is at least  $ \frac{(r+k)(r+1-k)}{2}$. 
Since $M$ has $(r+1)(r+1-k)$ columns, the system (\ref{hsisc}) has over the  quotient field of $ S = K[t,u] $ at most 
$$
(r+1)(r+1-k) - \frac{(r+k)(r+1-k)}{2} = \binom{r+2-k}{2} 
$$ 
independent solutions. Since we have already found that many independent solutions with entries in $S$ we have determined all solutions of the system (\ref{hsisc}). 

Using this information we get a lower bound for the dimension of $
H^0(C, {\mathcal N}_C)$ 

\begin{prop} \label{prop-normal-est} 
Let $ C \subset \pp^n$ be a rope with degree $ n-k > 2$ and genus $ g $ that is supported on a line.  Let $ \alpha_0, \dots,
\alpha_{r-k} $ be the degrees of the rows of the matrix $ A $
associated to $ C.$ Assume $\alpha_0 \geq 2$. 
Then we have 
$$ h^0(C, {\mathcal N}_C) \geq
(r+1)(2+k-g) - k^2 + \sum_{l, h = 0}^{r-k} \sum_{v=h}^{r-k}
\max(0, \alpha_l -\alpha_h -\alpha_v +2). 
$$
\end{prop}

\begin{proof} 
 We know that $ P^{ij}=0 $ for every $ i,j $ by Proposition
5.3, and that $ Q^s_l \in [S]_{\beta_s + \alpha_l}$ are arbitrary polynomials.  This gives  
$$ \sum_{s=1}^k
\sum_{l=0}^{r-k} ( \beta_s + \alpha_l +1) = (r+1)(k-g) - k^2  
$$ 
parameters. 

Now, we claim that different choices of $\lambda := (\lambda_0,\ldots,\lambda_r)^t$ led to different vectors ${\bf P} = B^t \lambda$. Indeed, if $B^t \lambda = 0$ then $\lambda$ is a linear combination of the columns of $A^t$. Since $\alpha_0 \geq 2$, the entries of $A$ have degree $\geq 2$. This gives a contradiction because $\lambda_0,\ldots,\lambda_r$ have degree one. 

Furthermore, since ${\bf Q}_j^l + a_{l j} \lambda = A^t {\bf c}^{l j}$, we see that for fixed $\lambda$ different choices of ${\bf c}^{l j}$ lead to different solutions of the system 3. Indeed, the matrix $A^t$ provides an injective homomorphism of $S$--modules. 

Summing up, we have shown that different choices of $(\lambda^t, {\bf c}^{0 0},\ldots,{\bf c}^{r, r-h})$ lead to different solutions of the system 3. Since $\lambda_0,\ldots,\lambda_r \in [S]_1$ can be chosen arbitrarily we get $2 (r+1)$ further parameters in $K$. 

The parameters $ c_h^{lj} $ have to 
satisfy the systems (\ref{sisc}).  Thus, they can be written as 
\begin{eqnarray*} (\dots, c_h^{l0},
\dots, c_h^{lr}, \dots) = \sum_{i=0}^{r-k} \sum_{i'=i}^{r- k}
d_{i,i'}^l (\dots, A_{i'}, \dots, A_i, \dots) + \hspace*{3cm} \\
\hfill  (\dots, 0,
\lambda_0, \dots, \lambda_r, 0, \dots) 
\end{eqnarray*} 
where $d_{i,i'}^l$ is any polynomial in $S$ of degree
$\alpha_l -\alpha_i -\alpha_{i'}+1.$  Hence, the $ c_h^{lj}$'s depend on 
$$ 
\sum_{l=0}^{r-k} \sum_{i=0}^{r-k} \sum_{i'=i}^{r-k} \max \{0,\; 
\alpha_l -\alpha_i -\alpha_{i'} +2\}
$$ 
parameters in $K$. 

Adding up the various contributions provides our claim. 
\end{proof}

The previous proof also shows. 

\begin{corol} \label{cor-normal} 
Adopting the notation and assumptions in Proposition \ref{prop-normal-est} assume in addition that,  for the rope $C$, Condition \ref{cond} is satisfied for every $\ffi \in H^0(C, {\mathcal N}_C)$. Then we have 
$$ h^0(C, {\mathcal N}_C) =
(r+1)(2+k-g) - k^2 + \sum_{l, h = 0}^{r-k} \sum_{v=h}^{r-k}
\max(0, \alpha_l -\alpha_h -\alpha_v +2). 
$$

\end{corol} 

\begin{proof} 
Condition \ref{cond} guarantees that we have found all the solutions of the systems 1.\ - 3.\ as shown in the discussion preceding Proposition \ref{prop-normal-est}.  
\end{proof} 

\section{Parameter spaces for ropes} \label{sec-para-spaces} 

The aim of this section is to describe some natural parameter
spaces for families of $ (n-k)$--ropes in $ \pp^n $
 of fixed genus $ g$ that are supported on a line and to study their
relations. This leads to results on the Hilbert scheme of ropes.

By Proposition \ref{prop-genus} we know that $ g = -\sum_{j=1}^k
\beta_j = - \sum_{i=0}^{r-k} \alpha_i.$ Moreover, a rope is uniquely determined 
when we fix the line that supports the rope and we know the
sequence (\ref{vectfield}) that defines the scheme structure of
the rope. We introduce some further notation. 

\begin{defi} 
Let  $ C \subset \pp^n$ be an $ (n-k)$--rope whose  sequence (\ref{vectfield}) is 
$$ 
0 \to
\oplus_{j=1}^k \oo_{\pp^1}(-\beta_j -1) \to \oo_{\pp^1}^{r+1}(-1) \to
\oplus_{i=0}^{r-k} \oo_{\pp^1}(\alpha_i-1) \to 0.  
$$
We may assume that 
$$
\alpha_0 \leq \alpha_1 \leq \ldots \leq \alpha_{r-k} \qquad \mbox{and}\qquad \beta_1 \leq \ldots \leq \beta_k. 
$$
In this case  we call $ {\bf \alpha} =(\alpha_0, \dots, \alpha_{r-k}) $ the  {\em right
type} and  $ {\bf \beta} = (\beta_1, \dots,
\beta_k) $ the  {\em left type} of the rope $C$. 
\end{defi}

Now, we can describe our first parameter space. 

\begin{prop} \label{prop-dim-Val} 
The $ (n-k)$--ropes with fixed right type $ {\bf \alpha}$  are
parameterized by an irreducible, smooth, quasi-projective scheme $ V_{\bf \alpha} $
of dimension 
\begin{equation} \label{dimalpha}
\dim V_{\bf \alpha} = (n-1)(n-k+1-g) - \sum_{i,j = 0}^{r-k}
{{\alpha_i - \alpha_j +1}\choose 1}. 
\end{equation} 
\end{prop} 

\begin{proof} 
According to \cite{NNS}, Proposition 2.9, two matrices $A, A'$ providing surjective morphisms 
$\oo_{\pp^1}^{r+1}(-1)
\stackrel{}{\lra} \oplus_{i=0}^{r-k} \oo_{\pp^1}(\al_i-1)$ of maximal rank lead via their syzygy matrices to the same rope if and only if $A = P \cdot A'$ where the matrix $P$ provides an automorphism of $\oplus_{i=0}^{r-k} \oo_{\pp^1}(\al_{i}-1)$. This is equivalent to the fact that $A^t$ and $(A')^t$ define the same subbundle of $\oo_{\pp^1}^{r+1}(1)$. Thus, denoting by $\Grass(1,n)$ the Grassmanian of lines in $\pp^n$ we get as parameter space 
$$ 
V_{\bf
\alpha} = \Grass(1,n) \times \left\{\mbox{subbundles of}\
\oo_{\pp^1}^{r+1}(1)\  \mbox{isomorphic to}\ \oplus_{i=0}^{r-k}
\oo_{\pp^1}(1 - \al_i) \right\}.  
$$ 
Hence, $V_{\bf \alpha}$ is smooth, irreducible and has dimension 
\begin{eqnarray*}
\dim V_{\bf \alpha} & = & 
2(r+1)+(r+1)\sum_{i=0}^{r-k}(\alpha_i+1) - \dim
\mbox{Aut}(\oplus_{i=0}^{r-k} S(1 - \alpha_i))  \\ 
&  = & 
(n-1)(n-k+1-g) - \sum_{i,j=0}^{r-k}{{\alpha_i-\alpha_j+1}\choose 1}
\end{eqnarray*}
because the supporting line fixes the base of $
S^{r+1}(-1)$. 
\end{proof} 

There is an analogous result using the
left type. We omit its proof. 

\begin{prop} \label{prop-Vbe} 
The $ (n-k)$--ropes with fixed left type $ {\bf \beta}$ are
parameterized by the irreducible, smooth, quasi-projective scheme 
$$ 
W_{\bf \beta} =
\Grass(1,n) \times \left\{ \mbox{subbundles of}\ \oo_{\pp^1}^{r+1}(1)\
\mbox{isomorphic to}\ \oplus_{j=1}^{k} \oo_{\pp^1}(1 - \beta_j) \right\}  
$$
of dimension 
$$ 
\dim W_{\bf \beta} = (n-1)(k+2-g) - \sum_{i,j=1}^{k} {{\beta_i - \beta_j + 1}\choose 1}.
$$ 
\end{prop}

\begin{obs} \rm Let $C \subset \pp^n$ be a non-degenerate $ (n-k)$--rope with left type $ {\bf \beta}$. Then its generic initial ideal
with respect to the degree reverse lexicographic order is 
$$
J=(x_0^2, x_0x_1, x_1^2, \dots, x_r^2, x_0t^{\beta_1}, \dots,
x_{k- 1}t^{\beta_k}) 
$$ 
because $ J \subset R$  is the only saturated,
Borel-fixed, monomial ideal compatible with the homogeneous ideal and the Hilbert function of $C$.  Hence, all the ropes parametrized by $W_{\bf \beta}$ correspond to points in the closure of the stratum of the Hilbert scheme $H_{n-k, g}(\pp^n)$ that parametrizes curves with initial ideal $J$ (cf.\ \cite{ns}).
\end{obs} 

Now, we begin to compare the parameter spaces. 

\begin{obs} \label{rem-trans}  
\rm It is not too difficult to see that the ropes parametrized by  $ V_{\bf \alpha}$  can have every possible left type   $ {\bf \beta} $ 
subject to the condition 
$$
| {\bf \beta}| :=  \sum_{j=1}^k \beta_j = -g = \sum_{i=0}^{r-k} \al_i.  
$$ 
In fact, this is a special case of \cite{schles}, Theorem 2.1.  Thus, every space $ W_{\bf \beta}$, where $ {\bf \beta} $ satisfies the condition above, has a point that corresponds to a point in $V_{\bf \al}$.  Similarly, varying a point  in $ W_{\bf \beta}$ we meet all the spaces $V_{\bf \al}$ such that the condition above is satisfied. 
Hence, the two types of parameter
spaces are transverse to each other, in some sense. 
\end{obs} 

The natural incidence relation between $
\cup_{\bf \alpha} V_{\bf \alpha} $ and $ \cup_{\bf \beta} W_{\bf
\beta}$ can be easily described. Let $(l, \cE)$ be a point of some $V_{\bf \alpha}$ and let $(l', \cE')$ be a point of some $W_{\bf \beta}$. The the incidence relation is just 
$$ 
{\mathcal H} = \left\{ ((l,\cE),(l',\cE'))
\s l = l',\; \cE' \cong \cE^* \right\}. 
$$
The condition $\cE' \cong \cE^* $  guarantees that there is an exact sequence 
$$ 
0 \to \oplus_{j=1}^k S(-\beta_j-1) \stackrel{\psi}\to
S^{r+1}(-1) \stackrel{\varphi}\to \oplus_{i=0}^{r-k} S(\alpha_i-1) 
$$ 
such that $\cE' = \widetilde{\im \psi}$ and $\cE^* = \widetilde{\im \ffi^*}$. 
\smallskip 

There is a third natural space parameterizing  $(n-k)$--ropes of genus $ g$.  It is the Hilbert scheme $\mbox{Hilb}_{n-k, g}(\pp^n) $ that parameterizes closed subschemes
in $ \pp^n $ with Hilbert polynomial $ p(z) = (n-k)z+1-g.$ Of
course, this Hilbert scheme contains not only ropes but many other
schemes with different geometric properties. It is natural
to investigate the relation among $ \cup_{\bf \alpha} V_{\bf
\alpha},$ $ \cup_{\bf \beta} W_{\bf \beta} $ and
$\mbox{Hilb}_{n-k, g}(\pp^n).$ Because of the universality of the
Hilbert scheme, there are  morphisms 
$$ j_{\bf \alpha} : V_{\bf
\alpha} \to \mbox{Hilb}_{n-k, g}(\pp^n)  \qquad \mbox{and}\qquad  j_{\bf \beta} :
W_{\bf \beta} \to \mbox{Hilb}_{n-k, g}(\pp^n) 
$$ 
that embed $ V_{\bf
\alpha} $ and $ W_{\bf \beta} $ into $\mbox{Hilb}_{n-k, g}(\pp^n)$ where $-g = |\al| = |\beta|$. 

Now, we prove a topological result about the images of these
morphisms. 

\begin{lemma} \label{lem-unique-comp} 
All $(n-k)$--ropes in $\pp^n$ with genus $g$ lie in the same irreducible component of $\mbox{Hilb}_{n-k, g}(\pp^n)$. More precisely, $ \bigcup_{{\bf |\alpha|} = -g} j_{\bf \alpha}(V_{\bf
\alpha}) $ is contained in one irreducible component of $\mbox{Hilb}_{n-k, g}(\pp^n).$
\end{lemma} 

\begin{proof} $ V_{\bf \alpha}$ being  irreducible, $ j_{\bf
\alpha}(V_{\bf \alpha}) $ is irreducible. Thus, it is enough to
show that for every two right types $ {\bf \alpha'} $
and ${\bf \alpha''}$ such that $ |{\bf \alpha'}| = |{\bf
\alpha''}| = -g,$  there is a flat family of ropes which
connects a point in $ j_{\bf \alpha'}(V_{\bf \alpha'}) $ and a
point in $ j_{\bf \alpha''}(V_{\bf \alpha''}).$ In order to construct the
family, we fix a line $ l \subset \pp^n$ and a left type $ {\bf \beta}$ such that $-g = |{\bf \beta}|$.  In $W_{\bf \beta}$ we choose two points $ (l,\widetilde{\im \psi'}) $ and $
(l,\widetilde{\im \psi''}) $ such that the morphisms $ \varphi',$ and $ \varphi''
$ corresponding to $ \psi',$ and $ \psi''$ via the incidence
relation above have type ${\bf \alpha'} $ and ${\bf \alpha''},$ i.e.\ $
(l,\widetilde{\im \varphi'}) \in V_{\bf \alpha'}$ and $ (l,\widetilde{\im \varphi''}) \in V_{\bf
\alpha''}.$ This is possible by Remark \ref{rem-trans}. 

Now, we define the family we are looking for. All the
ropes of the family are supported on the line $ l$ and the scheme
structure is defined via the left morphism of the sequence
(\ref{vectfield}). We set 
$$ 
{\mathcal F} =
\{ (l, s \psi' + (1-s) \psi'') | s \in U  \} 
$$
where $ U $ is the non-empty open subset of $ {\AAA}^1 $ defined by 
$U := \{s \in \AAA^1 \s 
\mbox{codim}(I_k(s \psi' + (1-s) \psi'')) = 2\}.$ Here, $I_k(\psi)$ denotes the determinantal ideal generated by the $k$-minors of the map $\psi$. Theorem \ref{charact} shows that ${\mathcal F}$ is a flat family and that the image
of ${\mathcal F}$ via $ j_{\bf \beta}(W_{\bf \beta}) $ in $
\mbox{Hilb}_{n-k, g}(\pp^n) $ connects $ j_{\bf \alpha'}(V_{\bf
\alpha'}) $ to $ j_{\bf \alpha''}(V_{\bf \alpha''}).$ 
\end{proof} 

The preceding result justifies. 

\begin{notation} \label{not-comp}
Let $ H_{n-k, g}(\pp^n) $ or shorter $ H_{n-k, g}$  be the irreducible
component of $ 
\mbox{Hilb}_{n-k, g}(\pp^n)$ that contains  
$j_{\bf \alpha}(V_{\bf \alpha}) $ for each ${\bf \alpha}$ with $-|\al| = g$. 
\end{notation} 

The next step is to find large families $V_{\al}$ and $W_{\beta}$ such that $j_{\al}(V_{\al}) \subset \mbox{Hilb}_{n-k, g}(\pp^n)$ and $j_{\beta}(W_{\beta}) \subset \mbox{Hilb}_{n-k, g}(\pp^n)$. We will be guided by  semicontinuity
of $ h^1(\pp^n, {\mathcal I}_C(z)) $ on $ H_{n-k, g}$. 
Recall that by Corollary \ref{raoC} the Rao function of a rope $C$ is 
$$ 
\rho_C(z) = \sum_{i=0}^{r-k}
{{z+\alpha_i}\choose 1} + \sum_{j=1}^k {{z-\beta_j}\choose 1}
-(r+1) {z \choose 1}.  
$$
It is determined by the left type and right type of $C$. Moreover, 
the left type $ {\bf \beta} $ does not contribute if $z \leq 0$, more precisely,   
$$ 
\rho_C(z) =
\sum_{i=0}^{r-k} {{z+\alpha_i}\choose 1} \qquad \mif z \leq 0, 
$$
while, for $ z>0,$ the right type $ {\bf \alpha}$ gives no
contribution because 
$$ \rho_C(z) = -k z -g + \sum_{j=1}^k
{{z-\beta_j}\choose 1} \qquad \mif z
> 0.
$$
 
 In order to specify the smallest Rao function of a rope in $H_{n-k, g}$ we introduce.

\begin{defi} \label{def-min-types}
Define integers $p_1, s_1$ by 
$$
-g = p_1 \cdot (r+1-k) + s_1 \qquad \mbox{where}\qquad 0 \leq s_1 \leq r-k. 
$$
Then we set $$
 {\bf \alpha}_{min} = (p_1,
\dots, p_1, p_1+1, \dots, p_1+1)
$$ 
where $ p_1+1 $ occurs $ s_1 $
times  while $ p_1 $ occurs $ r+1-k-s_1 $ times. 

Similarly, define integers $p_2, s_2$ by 
$$
-g = p_2 \cdot k + s_2 \qquad \mbox{where}\qquad 0 \leq s_2 < k. 
$$ 
We set 
$$ {\bf \beta}_{min} = (p_2, \dots, p_2, p_2+1, \dots,
p_2+1) 
$$ 
where $ p_2+1 $ occurs $ s_2 $ times, while $ p_2 $
occurs $(k-s_2)$ times. 
\end{defi} 

Obviously, ${\bf \alpha}_{min}$ is a right type  and ${\bf \beta}_{min}$ is a left type of a rope of genus $g$. They are rather particular. 

\begin{prop} \label{prop-small-coho} 
Let $\rho_{min}$ be the Rao function associated to the types $
{\bf \alpha}_{min} $ and $ {\bf \beta}_{min}$ (defined via the formula in Corollary \ref{raoC}). Then 
the Rao function of every $ (n-k)$--rope $ C \subset \pp^n$ of genus $ g $ 
satisfies 
$$
\rho_C(z) \geq \rho_{min}(z) \qquad \mbox{for all} \; z \in \ZZ. 
$$
\end{prop}

\begin{proof} Let $ {\bf \alpha'} = (\alpha_0, \dots, \alpha_{r-k}) $
and $ {\bf \alpha''} = (\alpha_0, \dots, \alpha_i+1, \dots,
\alpha_j-1, \dots,$ $ \alpha_{r-k}) $ be right types of some ropes such that $|{\bf \alpha'}| = |{\bf \alpha''}| = -g$. This means in particular that $\al_i + 2 \leq \al_j$. 
 We want to show that the Rao
function $ \rho_1 $ associated to $ {\bf \alpha'}$ is larger
than the Rao function $ \rho_2 $ associated to $ {\bf \alpha''}$ if $ z \leq 0.$ 

In fact, if $z \leq 0$ then we have 
$$ 
\rho_1 (z) = \sum_{h \not= i,j} {{z+\alpha_h}\choose 1}
+ {{z+\alpha_i}\choose 1} + {{z + \alpha_j}\choose 1} 
$$ 
and $$
\rho_2 (z) = \sum_{h \not= i,j} {{z+\alpha_h}\choose 1} +
{{z+\alpha_i+1}\choose 1} + {{z + \alpha_j - 1}\choose 1},
$$ 
thus
$$ \rho_1(z) - \rho_2(z) = {{z+\alpha_j-1}\choose 0} -
{{z+\alpha_i}\choose 0} \geq 0
$$ 
because $\al_i + 2 \leq \al_j$. 

Analogously, let $ {\bf \beta'} = (\beta_1, \dots, \beta_k) $ and
$ {\bf \beta''} = (\beta_1, \dots, \beta_i+1, \dots, \beta_j-1,
\dots, \beta_k) $ be left types such that $|{\bf \beta'}| = | {\bf \beta''}| = -g$. 
The Rao functions associated to these types are in case $z > 0$ 
$$ 
\rho_1(z) = -kz -g + \sum_{h \not= i,j}
{{z-\beta_h}\choose 1} + {{z- \beta_i}\choose 1} +
{{z-\beta_j}\choose 1} 
$$ 
and 
$$ 
\rho_2(z) = -kz -g + \sum_{h
\not= i,j} {{z-\beta_h}\choose 1} + {{z-\beta_i- 1}\choose 1} +
{{z-\beta_j+1}\choose 1}.$$ 
Their difference is 
$$ 
\rho_1(z) -
\rho_2(z) = {{z-\beta_i-1}\choose 0} - {{z- \beta_j}\choose 0}
\geq 0 \qquad \mif  z > 0.
$$
Hence, the Rao function is minimal  when the difference of two
entries of each type is as at most one. Our claim follows. 
\end{proof} 

Due to semicontinuity of the Rao function, we expect the parameter
space  $ V_{{\bf \alpha}_{min}}$ to have the largest dimension among the parameter spaces $V_{\al}$ with image in $H_{n-k, g}$. Now,  we compute its dimension. 

\begin{prop} \label{prop-dim-Valmin} 
$ \dim V_{{\bf \alpha}_{min}} = (n-1)(k+2-g) - k^2.$ 
\end{prop}

\begin{proof} 
According to the proof of Proposition \ref{prop-dim-Val} it suffices to compute 
 the dimension of the automorphism group of
$ \oplus S(\alpha_i-1) $ corresponding to the right type $
{\bf \alpha}_{min}.$ By Definition \ref{def-min-types} we have 
$$ 
\alpha_i = \left\{ \begin{array}{ccc} p_1 & \mbox{if } & i=0,
\dots, r-k-s_1 \\ p_1+1 & \mbox{if } & i= r+1-k-s_1, \dots, r-k
\end{array} \right., 
$$ 
thus $ \dim \mbox{Aut}(\oplus_{i=0}^{r-k}(\alpha_i-1)) = (r+1-k)^2.$
Hence, using formula (\ref{dimalpha}) an easy computation provides our claim. 
\end{proof} 

\begin{obs} \rm There is an analogous result for the left type. In fact, we get 
$$ 
\dim W_{{\bf \beta}_{min}} = (n-1)(k+2-g) -
k^2.
$$ \end{obs} 

Now, we are ready to study the component  $H_{n-k, g}(\pp^n) $ of
the Hilbert scheme $ \mbox{Hilb}_{n-k, g}(\pp^n) $ that contains 
ropes supported on lines (cf.\ Notation \ref{not-comp}). For the remainder of this paper we assume that the field $K$ is algebraically closed. 

First, we discuss double lines. Their maximal genus is $-1$ and the case of genus $-1$ is particular. 

\begin{obs} \label{rem-g--1} 
{\rm The component $H_{2,-1}(\pp^n)$ is generically smooth of dimension $4 (n-1)$. Its general element is a pair of skew lines. This is well-known (cf., e.g., \cite{Harris-Montreal}). }
\end{obs} 

In the general case, the characteristic of the ground field plays a role.  

\begin{thm} \label{thm-double-line} 
Let $g \leq -2,\; n \geq 3$ be integers. Then we have: 
\begin{itemize} 
\item[(a)] $\dim H_{2, g}(\pp^n_K) = (n-1) (3 - g) - 1$; 
\item[(b)] $H_{2, g}(\pp^n_K)$ is generically smooth if and only if either $\chr K \neq 2$ or $\chr K = 2$, $g = - 2$. 
\item[(c)] $H_{2, g}(\pp^n_K)$ is non-reduced if $\chr K = 2$ and $g \leq -3$.  
\end{itemize} 
Furthermore, a double line $C \subset \pp^n$  is obstructed if it has genus $\leq -3$ and $\chr K = 2$; otherwise unobstructed.  
\end{thm} 

\begin{proof} 
Every curve having degree two and genus at most $-2$ is a double line. Hence, Lemma \ref{lem-unique-comp} shows that the dimension of $H_{2, g}(\pp^n)$ is the maximum of $\dim V_{\al}$ where $\al$ is the right type of a double line. Proposition \ref{prop-dim-Val} and Proposition \ref{prop-dim-Valmin} provide the claim about the dimension. 

Corollary \ref{h0nc2} gives the dimension of the tangent space of  $H_{2, g}(\pp^n_K)$ at every point corresponding to a double line $C$. But the proof of (a) shows that the general curves of $H_{2, g}(\pp^n)$ are the ropes in $V_{\al_{min}}$. Thus, $H_{2, g}(\pp^n)$ is generically smooth if and only if every $C \in V_{\al_{min}}$ satisfies $(n-1) (3 - g) - 1 = h^0(C, {\mathcal N}_C)$. Hence, Corollary \ref{h0nc2} implies claim (b) as well as (c) because the curves form an open subset of $Hilb_{2, g}(\pp^n)$. 
\end{proof} 

%Of course, $ \dim V_{{\bf \alpha}_{min}} $ is a lower bound for the
%dimension of $ H_{n-k, g}$ An upper bound for the dimension of $ H $ is $
%h^0(C,{\mathcal N}_C) $ where $C$ is an $ (n-k)$--rope $ C $ of genus $ g = - |\al_{min}|.$ 

It remains to consider ropes whose degree is at least three. In this case our main result about the component $ H_{n-k, g}(\pp^n)$ is. 

\begin{thm} \label{thm-main-ropes} 
Let $d, g, n$ be integers such that $3 \leq d \leq n-1$ and $g \leq d-n$. 
If either $ g \leq -3(d-1) $ or $ -g = 2(d-1)$ then the component $ H_{d, g}(\pp^n)$ is
generically smooth of dimension $ (n-1)(d+1-g) - (d-1)^2$  and its 
general element is a non-degenerate rope. 
\end{thm} 

\begin{proof} 
Put $n - k := d$. 
Of course, $ \dim V_{{\bf \alpha}_{min}} = (n-1)(k+2-g) - k^2$ is a lower bound for the 
dimension of $ H_{n-k, g}(\pp^n)$. Hence, it is enough to prove that the tangent space of $ H_{n-k, g}(\pp^n) $ at a suitable point  has dimension $ (n-1)(k+2-g) - k^2,$ i.e.\  
$$ 
h^0(C,{\mathcal N}_C) = (n-1)(k+2-g) - k^2 
$$ 
for some $(n-k)$--rope $ C \subset \pp^n$ of genus $ g.$ 

To this end we consider ropes of varying embedding dimension. Let $C \subset \pp^n$ be a rope having the right type $\al = (\al_1,\ldots,\al_{s+w}) := (p,\ldots,p, p+1,\ldots,p+1)$ with $w \geq 1$ occurrences of $p \geq 1$ and $s \geq 0$ occurrences of $p+1$. Then (cf.\ Section \ref{sec-char}) $C$ has degree $s + w + 1$ and genus $g := -[(s+w) p + s]$. Moreover, we must have $s+w < r+1$. 
The rope $C$ is non-degenerate if and only if its left type $\beta = (\beta_1,\ldots,\beta_k)$ satisfies $\beta_1 \geq 1$. In this case we get $-g = \beta_1 + \ldots + \beta_k \geq k = r+1 - (s+w)$ which implies 
$$
s+w \leq r \leq (s+w) \cdot (p+1) + s - 1. 
$$

We now fix the right type $\alpha$ but let the linear span of the ropes vary. Our goal is to show.  

{\it Claim}: For every $r$ in the range above  there is a
non-degenerate rope $C \subset \pp^{r+2}$ with right type $\al$ such
that Condition \ref{cond} is satisfied for every $\ffi \in H^0(C,
{\mathcal N}_C)$. (We use the notation introduced at the beginning of
Section \ref{sec-normal-sheaf}. Recall in particular that $P^s$
denotes just a polynomial with index $s$. It is not the $s$-th power
of $P$.)  

We begin by specifying a suitable rope in $\pp^{r+2}$.  According to  Theorem \ref{charact} it suffices to describe an $(s+w) \times (r+1)$ matrix $A_{r}$. We do this with the help of three types of matrices. 

First, we define the $(s+w) \times (s+w+1)$ matrix $\tilde{A}_0  = (a_{ij})$ by 
$$
a_{ij} = \left\{ \begin{array}{ccl} 
t^p & \  & \mif i = j,\; i=0, \dots,w-1; \\
t^{p+1} & \  & \mif i = j,\; i=w, \dots, s+w-1; \\ 
u^p & \  & \mif j=i+1,\; i=0, \dots,w-1; \\ 
u^{p+1} & \  & \mif j=i+1,\; i=w, \dots, s+w-1; \\ 
0 & \  & \mbox{otherwise}
\end{array} \right. 
$$ 
If $0 \leq j \leq s+w$ we denote by $\tilde{A}_{j}$ the matrix that is obtained from $\tilde{A}_0$ by deleting its last $j$ rows. 

Second, we put 
$$
A_i' := (t^{p+1}, t^p u,\ldots,t^{p+2-i} u^{i-1}, u^{p+1}) \qquad (1 \leq i \leq p+1) 
$$
and 
$$
A_i'' := (t^{p}, t^{p-1} u,\ldots,t^{p+1-i} u^{i-1}, u^{p}) \qquad (1
\leq i \leq p).   
$$
Note that the format of these matrices is $1 \times (i+1)$. 

Now, we define the matrices $A_r$. If $r = s+w$ then we set $A_r := \tilde{A}_0$. 

Let $s + w < r \leq s+w + s (p+1)$. Define integers $i, j$ by 
$$
r  := s+w + j (p+1) + i \qquad \mbox{and} \qquad 1 \leq i \leq p+1. 
$$ 
Then we put
$$
A_r := \left ( \begin{array}{ccccc} 
\tilde{A}_{j+1} \\ 
& A_i' \\
& & A_{p+1}' \\
& & & \ddots \\
& & & & A_{p+1}' 
\end{array} \right )
$$
where all non-specified entries are zero. Observe that the block $A_{p+1}'$ occurs $j$ times. 

If $s (p+2) + w < r \leq (s+w) (p+1) + s - 1$ then we define integers $i, j$ by 
$$
r  := s (p+2)+w + j p + i \qquad \mbox{and} \qquad 1 \leq i \leq p. 
$$ 
Then we set 
$$
A_r := \left ( \begin{array}{cccccccc} 
\tilde{A}_{j+s+1} \\ 
& A_i'' \\
& & A_p''\\
& & & \ddots \\
& & & & A_p''\\
& & & & & A_{p+1}' \\
& & & & & & \ddots \\
& & & & & & & A_{p+1}' 
\end{array} \right )
$$
where block $A_{p+1}'$ occurs $s$ times and the block $A_p''$ appears $j$ times. 

Furthermore, we define $B_r$ as a syzygy matrix of $A_r$. It has also a block structure. 
It is easy to see that the maximal minors of $A_r$ and $B_r$, respectively, generate an ideal of codimension two. Hence, the ideal $ ((I_L)^2, [x_0, \dots, x_r] B_r)$ defines a non-degenerate rope $C \subset \pp^{r+2}$ with the prescribed right type $\al$. It remains to check that  Condition \ref{cond} is satisfied. 

To this end we induct on $r \geq s+w$. Let $r := s+w$. Then the entries of $B_r$ are the maximal minors of $ A_{e},$ i.e.\ 
$$ 
B_{e} = \left (u^{-g}, -t^pu^{-g-p}, t^{2p}u^{-g-2p}, \dots,
(-1)^{r-1}t^{-g-p- 1}u^{p+1}, (-1)^r t^{-g} \right)^t.
$$ 
Furthermore,  in this case the vector $ {\bf P} $ has just one entry. We denote it by $P$ to simplify notation a bit. Then for $C$, the equations of the linear system $ 3.$ in Proposition \ref{prop-lin-systems} with $ l=0$ are
\begin{equation} 
\label{ad} \left\{ \begin{array}{lcc} u^{-g} Q^{00}_0 +
\dots + (-1)^r t^{-g} Q^{0r}_0 + t^p P & = & 0 \\ u^{-g} Q^{01}_0
+ \dots + (-1)^r t^{- g} Q^{1r}_0 + u^p P & = & 0 \\ u^{-g}
Q^{02}_0 + \dots + (-1)^r t^{-g} Q^{2r}_0 & = & 0 \\ \vdots \\
u^{-g} Q^{0r}_0 + \dots + (-1)^r t^{-g} Q^{rr}_0 & = & 0
\end{array} \right. 
\end{equation}
The first equation shows that $ t^p $ divides $ Q^{00}_0 $, thus $Q^{00}_0 = t^p q_0$ for some $q_0 \in S$.  Analogously, the $ (i+1)$--st equation provides for $ i=2, \dots, r$  that $Q^{0i}_0 = t^p q_i $ for some $q_i \in S$. 
If we multiply the first equation by $ u^p,$ the second one by $
-t^p $ and then add the two equations, we get 
\begin{eqnarray*}  
u^{-g+p}t^p q_0 -
t^pu^{-g}Q^{01}_0 + t^{3p}u^{-g-p}q_2 + \dots + (-1)^rt^{-
g+p}u^p q_r \hspace*{3cm} \\ 
\hfill - t^pu^{-g}Q^{01}_0 + t^{2p}u^{-g-p}Q^{11}_0 +
\dots + (-1)^{r+1}t^{- g+p} Q^{1r}_0 = 0. 
\end{eqnarray*} 
Assume now that $\chr K \neq 2$. Then the last equation shows that 
$ t^p $ divides $ u^{-g+p}q_0 - 2 u^{-g} Q^{01}_0,$ i.e.\ $
u^{-g+p}q^0 - 2 u^{-g} Q^{01}_0 = t^p \tau'$ for some $\tau \in S$. 
It follows that $ u^{-g} $
divides $ \tau'$. Hence, there is a $\tau \in S$ such that $ Q^{01}_0 = \frac12 (u^p q_0
+ t^p \tau).$ Plugging this into the first equation of (\ref{ad})
we obtain 
$$ 
\frac12 u^{-g}t^p q_0 - \frac12 t^{2p}u^{-g-p} \tau +
t^{3p}u^{-g-2p} q_2 + \dots + (-1)^rt^{-g+p} q_r + t^p P = 0.  
$$ 
After cancelling  $ t^p $ we see that $P$ is a linear combination of the rows of $B_e$, i.e.\ Condition \ref{cond} is satisfied. 

If $\chr K = 2$ a similar argument shows the claim. 
This completes the initial step of the induction. 

Let $r > s+w$. We carry out  the induction step only in case $r = s+w+1$.  It shows the ideas. The general case is only notationally more complicated. 

Thus, let $r := s+w+1$. Then the syzygy matrix of $A_r$ is 
$$ 
B_{r} = \left( \begin{array}{cc} \widetilde{(B_{r-1})}_{1} & 0 \\ 0
& B'
\end{array} \right) 
$$
where $ \widetilde{(B_{r-1})}_{1}$ is the column obtained from $B_{r-1} $ by deleting the last row and then cancelling the common factor $u^{p+1}$ of  the entries,  and 
$$ 
B' = \left \{ \begin{array}{cl} 
(u^{p+1}, -t^{p+1})^t & \mif s > 0 \\
(u^p, -t^p)^t & \mif s = 0. 
\end{array} \right.
$$
For the rope defined by $ ((I_L)^2, [x_0, \dots, x_r] B_{r})$ the vector ${\bf P}$ is $(P^1, P^2)$ and the linear system $ 3.$ becomes 
$$
\left\{
\begin{array}{l} 
(\widetilde{(B_{r-1})}_{1})^t {\bf Q}^j_l + a_{lj} P^1 = 0 
\\ (B')^t {\bf Q}^j_l + a_{lj} P^2 = 0 \end{array} \right.
$$ 
Using the first set of equations we see, as in the previous case $r = s+w$ above, that $ P^1 $ belongs to the ideal
generated by the first $r$ entries of the first column of $B_r$. To conclude for $P^2$ we use the part of the system where $l=0$. It reads as 
\begin{eqnarray} \label{eq-ind} 
B^t_r {\bf Q}_0^j + a_{0 j} \cdot \left (\begin{array}{c}
P^1 \\
P^2
\end{array} \right ) = 0. 
\end{eqnarray}
By the definition of $A_r$, this provides for $j = 2,\ldots,r$ 
$$
B^t_r {\bf Q}_0^j = 0, 
$$ 
hence ${\bf Q}_0^j = A_r^t \cdot {\bf d}_j$ for some vector ${\bf d}_j = (d_{0 j},\ldots,d_{r-2, j})$, in particular 
$$
Q^{0 j}_0 = t^p d_{0 j} \qquad \fall j = 2,\ldots,r. 
$$
For $j = 0,\; l = 0$ the system (\ref{eq-ind}) gives the subsystem  
$$
B^t_r {\bf Q}_0^0 + t^p \cdot \left (\begin{array}{c}
P^1 \\
P^2
\end{array} \right ) = 0.
$$
Using $Q_0^{0 j} = Q_0^{j 0}$ the second equation provides 
$$
u^{p+1} t^p d_{0, r} - t^{p+1} t^p d_{0, r+1} +  t^p P^2 = 0. 
$$
Cancelling $t^p$ the claim follows for $P^2$. (Using the first instead of the second equation would give an alternative, more explicit proof of the claim for $P^1$.) This completes the proof of the claim. 
\smallskip 

Let now $d, g, n$ be integers as in the statement. 
Then the claim shows that  the component $H_{n-k, g}(\pp^n)$ contains a non-degenerate rope $C$ to which Corollary \ref{cor-normal} applies because the assumption $-g \leq 2 (n-1-k)$ implies $p = \al_0 \geq 2$. Thus, we get 
$$ 
h^0(C, {\mathcal N}_C) = (r+1)(2+k-g) - k^2 + \sum_{l,h}^{0,\dots,
r-k} \sum_{v=h}^{r-k} \max(0, \alpha_l -\alpha_h -\alpha_v +2) -
\delta_{1,\alpha_0}.  
$$ 
It follows that $h^0(C, {\mathcal N}_C) = (r+1)(2+k-g) - k^2$ if and
only if $\al_{r-k} - 2 \al_0 \leq -2$. The latter is true if $\al =
(2,\ldots,2)$, i.e.\ $-g = 2 (n-1-k)$, or $\al_0 \geq 3$, i.e.\ $-g
\leq 3 (n-1-k)$.  
Comparing with Proposition \ref{prop-dim-Val} we see that the rope $C$
is unobstructed.  
\end{proof} 

\begin{obs} \rm Theorem \ref{thm-main-ropes} shows that there is an open subset  in $H_{d, g}(\pp^n)$ consisting of ropes such that the global sections of the normal sheaf computed in Section \ref{sec-normal-sheaf} are all the global sections. In other words, Condition \ref{cond} is an open condition in $H_{d, g}(\pp^n)$ (cf.\ Corollary \ref{cor-normal}).  
\end{obs}

Similarly, our knowledge of the normal sheaf allows to characterize singular points of $H_{d, g}(\pp^n)$. 

\begin{corol} 
Let $d, g, n$ be as in Theorem \ref{thm-main-ropes}. Then a
rope $ C $ in $H_{d, g}(\pp^n)$ is obstructed if and only if either
the corresponding 
system $ 3. $ has more solutions than the ones we computed (i.e.\ the
inequality in Proposition \ref{prop-normal-est} is strict), or
$\alpha_{r-k} - 2 \alpha_0 \geq - 1.$  
\end{corol}

Recall that we mean by a curve always a locally Cohen-Macaulay $1$-dimensional scheme. Thus, if $g \leq -2$ then the only curves in $Hilb_{2, g}(\pp^n)$ are $2$--ropes. If $3 \leq d \leq n-1$ then $Hilb_{d, g}(\pp^n)$ does not contain only $3$--ropes even if $g \ll 0$. But if we restrict to the component $H_{d, g}(\pp^n)$ we get an analogous result. 

\begin{prop} \label{prop-only-ropes} 
Let $d, g, n$ be integers such that $3 \leq d \leq n-1,\; g \leq d-n$, and 
either $ g \leq -3(d-1) $ or $ -g = 2(d-1)$. Then all curves in $H_{d, g}(\pp^n)$ are ropes. 
\end{prop} 

\begin{proof}
According to Theorem \ref{thm-main-ropes}, the general element of
$H_{d, g}(\pp^n)$ is a rope. Hence, every curve in the closure of the
open subset of ropes is the limit of a flat, $1$-dimensional family
$\cF$ of ropes. Invoking projective transformations we may assume that
all the ropes in $\cF$ are supported on the {\em same} line $L$. But
then the limit $C$ of $\cF$ satisfies $I_C \supseteq (I_L)^2$, thus
$C$ is a rope.  
\end{proof} 

There are, of course, $1$-dimensional schemes in $H_{d, g}(\pp^n)$
that are not locally Cohen-Macaulay. It suffices to consider a rope on
a line $L$ and to deform its associated matrix $B$ to one, say $B'$,
that drops rank at some point of $L$. Then the ideal $((I_L)^2,
[x_0,\ldots,x_r] B')$ defines a non-locally Cohen-Macaulay scheme in
$H_{d, g}(\pp^n)$.  

\begin{obs} \rm
Let $\cF$ be the family of $1$-dimensional subschemes that are the
union of a double line of genus $g$ in $\pp^n$ and a point. Note that
its elements have genus $g-1$. We have just seen that $\cF$ and $H_{2,
  g-1}(\pp^n)$ meet. The dimension of $\cF$ is 
$$
\dim \cF = \dim (H_{2, g}(\pp^n)) + n = \dim (H_{2, g-1}(\pp^n)) +1. 
$$ 
It follows that the general element of $\cF$ is not contained in
$H_{2, g-1}(\pp^n)$. 
\end{obs} 

%%%%%%%%%%%%%%%%%%%%%%%%%%%%%%%%%%%%%%%%%%%%%%%%%%%%%%%%%%%%%%%%%%%%%%%%%%%%%

\end{document}